\documentclass[titlepage,final,11pt]{article}
\usepackage{theorem}
\usepackage{enumerate}
\usepackage{array}
\usepackage[reqno]{amsmath}
\usepackage{amssymb}
\usepackage{mathtools}
\usepackage{latexsym}
\usepackage{makeidx}
\usepackage{fancybox}
\input epsf.sty
\usepackage[dvips]{graphicx,color}
\usepackage{showlabels}
\usepackage{subcaption}
\usepackage{hyperref}

\theoremstyle{change}
\newtheorem{proclaim}{PROCLAIM}[section]

\numberwithin{equation}{section}

\outer\def\proclaim #1. #2\par{\medbreak \noindent{\bf#1.\enspace}{\sl#2}\par
  \ifdim\lastskip<\medskipamount
  \removelastskip\penalty55\medskip\fi}
\def\state #1. { \noindent{\bf#1.\enspace}}
\def\algo #1. { \noindent{\bf#1.\enspace}}

\newcommand{\comp}{\,{\raise 1pt \hbox{$\scriptstyle\circ$}}\,}

\newcommand{\natnums}{{{\rm l} \kern -.13em {\rm N} }}

\newcommand{\snats}{{I\kern -.29em N}}
\newcommand{\rats}{{Q\kern -.64em \raise 1pt \hbox{$\scriptstyle |$}\;\,}}
\newcommand{\srats}
	{{Q\kern -.56em \raise 1.2pt \hbox{$\scriptscriptstyle /$}\,}}
\newcommand{\ints}{Z\kern -.46em Z}

\newcommand{\pluss}{\hskip1pt \raise1pt\vbox{\hrule width6pt \vskip1pt \hrule
                    width6pt} \kern-4pt{\lower1pt\hbox{\vrule height6pt
		    \kern1pt\vrule height6pt}}\hskip5pt}
\newcommand{\eop}
	{\hfill{$\vcenter{\hrule height1pt \hbox{\vrule width1pt height5pt
   	 \kern5pt \vrule width1pt} \hrule height1pt}$} \medskip}

\newcommand{\setd}{{ d \kern -.15em l}}
\newcommand{\hatsetd}{ d \hat{\kern -.15em l }}

\renewcommand{\epsilon}{\varepsilon}
\renewcommand{\phi}{\varphi}

\newcommand{\tto}{\;{\lower 1pt \hbox{$\rightarrow$}}\kern -12pt
           \hbox{\raise 2.5pt \hbox{$\rightarrow$}}\;}
\newcommand{\overto}[1]{\,{\raise 0pt\hbox{$\rightarrow$}}\kern -9pt
     \hbox{\lower 3pt \hbox{$\scriptscriptstyle#1$}}\hskip6pt}
\newcommand{\underto}[1]{\,{\lower 1pt\hbox{$\rightarrow$}}\kern -9pt
     \hbox{\raise 4pt \hbox{$\,\scriptscriptstyle#1$}}\hskip7pt}
\newcommand{\bigoverto}[1]{{\raise 0pt\hbox{$\,\longrightarrow$}}\kern -16pt
     \hbox{\lower 3pt \hbox{$\scriptscriptstyle#1$}}\hskip4pt}
\newcommand{\bigunderto}[1]{\,{\lower 1pt\hbox{$\longrightarrow$}}\kern -16pt
     \hbox{\raise 4pt \hbox{$\,\scriptscriptstyle#1$}}\hskip6pt}
\newcommand{\bigbigto}[2]{\,{\raise 0pt\hbox{$\,\longrightarrow$}}\kern -16pt
     \hbox{\lower 3pt \hbox{$\scriptscriptstyle#2$}}\kern -10pt
     \hbox{\raise 4pt \hbox{$\,\scriptscriptstyle#1$}}\hskip7pt}
\newcommand{\downto}{{\raise 1pt \hbox{$\scriptscriptstyle \,\searrow\,$}}}
\newcommand{\upto}{{\raise 1pt \hbox{$\scriptscriptstyle \,\nearrow\,$}}}

\newcommand{\notimply}
	{\quad\hbox{$\Longrightarrow \kern -14pt {/}$}\hskip6pt\quad}

\newcommand{\lto}{\,{\lower 1pt\hbox{$\rightarrow$}}\kern -10pt
     \hbox{\raise 4pt \hbox{$\, \scriptstyle l$}}\hskip7pt}
\newcommand{\eto}{\,{\lower 1pt\hbox{$\rightarrow$}}\kern -11pt
     \hbox{\raise 4pt \hbox{$\, \scriptstyle e$}}\hskip7pt}
\newcommand{\hto}{\,{\lower 1pt\hbox{$\rightarrow$}}\kern -11pt
     \hbox{\raise 4pt \hbox{$\, \scriptstyle h$}}\hskip7pt}
\newcommand{\pto}{\,{\lower 1pt\hbox{$\rightarrow$}}\kern -11pt
     \hbox{\raise 4.5pt \hbox{$\, \scriptstyle p$}}\hskip7pt}
\newcommand{\cto}{\,{\lower 1pt\hbox{$\rightarrow$}}\kern -11pt
     \hbox{\raise 4pt \hbox{$\, \scriptstyle c$}}\hskip7pt}
\newcommand{\gto}{\,{\lower 1pt\hbox{$\rightarrow$}}\kern -11pt
     \hbox{\raise 4.5pt \hbox{$\, \scriptstyle g$}}\hskip7pt}
\newcommand{\sto}{\,{\lower 1pt\hbox{$\rightarrow$}}\kern -11pt
     \hbox{\raise 4pt \hbox{$\, \scriptstyle s$}}\hskip7pt}
\newcommand{\awto}{\,{\lower 1pt\hbox{$\rightarrow$}}\kern -15pt
     \hbox{\raise 4pt \hbox{$\, \scriptstyle aw$}}\hskip7pt}
\def\Nto{\,{\raise 1pt\hbox{$\rightarrow$}}\kern -13pt
     \hbox{\lower 3pt \hbox{$\, \scriptstyle N$}}\hskip7pt}
\def\Cto{\,{\raise 1pt\hbox{$\rightarrow$}}\kern -14pt
     \hbox{\lower 3pt \hbox{$\, \scriptstyle C$}}\hskip7pt}
\def\fto{\,{\raise 1pt\hbox{$\rightarrow$}}\kern -14pt
     \hbox{\lower 3pt \hbox{$\, \scriptstyle f$}}\hskip7pt}

\newcommand{\low}[1]{{\lower1pt \hbox{$\scriptstyle #1$}}}
\newcommand{\loww}[1]{{\lower2pt \hbox{$\scriptstyle #1$}}}
\newcommand{\high}[1]{{\raise1pt \hbox{$\scriptstyle #1$}}}

\newcommand{\cA}{{\cal A}}
\newcommand{\cB}{{\cal B}}

\newcommand{\cD}{{\cal D}}

\newcommand{\cF}{{\cal F}}

\newcommand{\cJ}{{\cal J}}

\newcommand{\cL}{{\cal L}}

\newcommand{\cR}{{\cal R}}
\newcommand{\cS}{{\cal S}}
\newcommand{\cT}{{\cal T}}

\newcommand{\cV}{{\cal V}}
\newcommand{\cW}{{\cal W}}

\newcommand{\cY}{{\cal Y}}

\newcommand{\nnmin}{\mathop{\rm minimize}}

\newcommand{\lwdy}[2]{\mathrel{\mathop
        {\raisebox{0.1ex}{\null$#1$}}{\hbox{\kern -1.0em
	{\raisebox{-0.8ex}{$\scriptstyle{\;\to #2}$}}}}}}
\newcommand{\lwwdy}[2]{\mathrel{\mathop
        {\raisebox{0.2ex}{\null$#1$}}{\hbox{\kern -1.0em
	{\raisebox{-1.1ex}{$\scriptstyle{\;\to #2}$}}}}}}
\newcommand{\slwdy}[2]{\scriptsize{{\mathrel{\mathop
        {\raisebox{0.1ex}{\null$#1$}}{\hbox{\kern -1.0em
	{\raisebox{-0.8ex}{$\scriptstyle{\;\to #2}$}}}}}}}}
\newcommand{\slwwdy}[2]{\scriptsize{{\mathrel{\mathop
        {\raisebox{0.2ex}{\null$#1$}}{\hbox{\kern -1.0em
	{\raisebox{-1.1ex}{$\scriptstyle{\;\to #2}$}}}}}}}}

\definecolor{lightgray}{gray}{0.75}
\definecolor{myred}{rgb}{0.55,0,0}
\definecolor{myblue}{rgb}{0,0,0.5} 
\definecolor{mygreen}{rgb}{0,0.5,0} 
\definecolor{purple}{rgb}{0.5,0,0.5} 
\definecolor{turq}{rgb}{0,0.805,0.816} 
\definecolor{maroon}{rgb}{0.51,0,0}
\definecolor{MAROON}{rgb}{0.51,0,0}
\definecolor{redor}{rgb}{0.78,0.078,0.078}
\definecolor{dgreen}{rgb}{0,0.3,0}

\newcommand{\bcdot}{\,{\raise .2ex \hbox{$\centerdot$}}\,}

\usepackage{graphics}
\usepackage{amsmath}
\usepackage{xcolor}

\usepackage[title]{appendix}
\usepackage{epsfig}
\usepackage{array}
\newcolumntype{P}[1]{>{\centering\arraybackslash}p{#1}} 
\usepackage{caption}
\usepackage{float}
\usepackage{mathtools}
\usepackage{comment}
\usepackage{url}

\newcommand{\ML}[1]{{\color{blue}{#1}}}

\begin{document}


\begin{center}
\begin{large}
{\bf Multi-Agent Search for a Moving and Camouflaging Target}
\smallskip
\end{large}
\vglue 0.7truecm
\begin{tabular}{ccc}
  {\sl Miguel Lejeune }  &  {\sl Johannes O. Royset} &  {\sl Wenbo Ma}\\
  School of Business & Operations Research Department & School of Business\\
  George Washington University & Naval Postgraduate School & George Washington University\\
  mlejeune@gwu.edu & joroyset@nps.edu & wenboma2011@gmail.com
\end{tabular}

\vskip 0.2truecm

\end{center}

\vskip 1.3truecm

\noindent {\bf Abstract}. \quad In multi-agent search planning for a randomly moving and camouflaging target, we examine heterogeneous searchers that differ in terms of their endurance level, travel speed, and detection ability. This leads to a convex mixed-integer nonlinear program, which we reformulate using three linearization techniques. We develop preprocessing steps, outer approximations via lazy constraints, and bundle-based cutting plane methods to address large-scale instances. Further specializations emerge when the target moves according to a Markov chain. We carry out an extensive numerical study to show the computational efficiency of our methods and to derive insights regarding which approach should be favored for which type of problem instance.

\vskip 0.5truecm

\halign{&\vtop{\parindent=0pt
   \hangindent2.5em\strut#\strut}\cr
                         \cr\cr

\cr }

\baselineskip=15pt

\vspace{-0.7in}

\section{Introduction}\label{sec:1}

Search for a randomly moving target in a discrete environment is challenging because the probability for detecting the target during a look at a particular location depends on the time of the look and the allocation of earlier looks. Thus, the optimization of searcher paths through discrete time and space results in difficult nonlinear problems with integer variables. Operational constraints on the searchers related to travel speed, endurance, and deconfliction further complicate the problem. In this paper, we formulate a mixed-integer nonlinear program (MINLP) that accounts for these factors. Given a planning horizon, it prescribes an optimal path for each searcher that maximizes the probability of detecting a randomly moving target that might camouflage, or not, and thus is even less predictable. We present a new linearized model and extend two others to account for operational constraints and heterogenous searchers. In an effort to reduce computing times, we develop a preprocessing technique, implement a lazy-constraint scheme within an outer-approximation solution method, and construct three cutting plane algorithms. Extensive numerical simulations demonstrate some of the modeling possibilities and indicate the most effective computational strategies in various settings.

Problems of the kind modeled in this paper arise in search-and-detection operations (see \cite{AbizeidMorinNilo.19} and \cite[Chapter 7]{Washburn.02} for a discussion of tools used by the U.S Coast Guard and the U.S. Navy), in counter-drug interdiction \cite{PietzRoyset.13,PietzRoyset.15,Zhangetal.20}, and in counter-piracy operations \cite{bourque2019}. It is also increasingly likely that planners in the near future will need algorithms for guiding large groups of autonomous systems as they carry out various search tasks, for example in underground environments \cite{DarpaUnder}.

The literature on search problems is extensive; see the reviews \cite{Ding.18,raap2019} as well as the monographs \cite{Washburn.02,Stone.04,StoneRoysetWashburn.16}. We assume a randomly moving target and {\em not} one that reacts or adapts to the searchers as seen, for example, in \cite{WashburnWood.95,Pfeiff.09} and \cite[Chapter 7]{StoneRoysetWashburn.16}. Thus, we broadly face the problem of optimizing a parameterized Markov decision process \cite{DimitrovMorton.09}, but can still avoid the formulation of a dynamic program and associated computational intractability as long as false positive detections are not considered. This fact is well-known and, at least, can be traced back to  \cite{Stewart.79}.

Specialized branch-and-bound algorithms using expected number of detections in bound calculations \cite{Washburn.98,LauHuangDissanayake.07,SatoRoyset.09} are effective when optimizing for a single searcher. Recently, this has been extended to multiple homogeneous searchers using minimum-cost flow computations to generate bounds \cite{bourque2019}. In the case of multiple searchers, cutting planes (constructed using either tangent or secant lines) furnish linear approximations that can be refined adaptively and lead to exact algorithms \cite{Royset2010}. The computational cost of identifying cuts tends to be significant if the target path can be any one of a large number of possible paths. This is reduced significantly when the target paths are governed by a Markov chain due to convenient formulas developed in \cite{Brown.80}; see also \cite{Royset2010}. Recent efforts toward developing cutting plane methods include \cite{delavernhe2021g}, but there exactness is sacrificed to achieve shorter computing times. The resulting algorithm uses a greedy heuristic to build the linear approximations.

A cutting plane approach can be viewed as a linearization of the actual problem when ``all'' cuts are included in the master problem from the outset. At least conceptually, this produces a direct solution approach: solve the master problem with all cuts included as proposed in \cite{Royset2010}. When the target moves according to a Markov chain, then one can also achieve another linearization through direct modeling of the evolution of the (posterior) probability of having the target in a particular location \cite{Royset2010}. This linearization approach is refined in \cite{berger2021near} under the assumption that the travel times between locations are always one time period and the searchers are homogeneous. This effort includes path splitting mitigation strategies for the continuous relaxation of the resulting mixed-integer model, variable elimination by switching to a focus on the terminal time period in the objective function, and implementation of a receding horizon strategy.

The literature also includes branch-and-bound algorithms that solve sequences of convex subproblems \cite{EagleYee.90} and many heuristics \cite{DellEagleMartinsSantos.96,Grundel.05,WongBourgaultFurukawa.05,RiehlCollinsHespanha.07,HollingerSingh.08,lanillos2012,AbizeidMorinNilo.19}, but they lack optimality guarantees. Routing of constrained searchers in discrete time and space has similarities with (team) orienteering and related reward-collecting vehicle routing problems; see, e.g., \cite{RoysetReber.09,PietzRoyset.13,ChoBatta.21,MoskalDasdemirBatta.23}. These problems often emphasize operational constraints such as time-windows for accomplishing tasks, limits on endurance and capacity, and deconflication among multiple agents.

In this paper, we also include operational constraints about endurance and deconflication, and hint to other possibilities that can be added with relative ease. In contrast to \cite{berger2021near}, which numerically examines one and two searchers, we study up to 50 searchers. We also allow for different types of searchers; their sensors, endurance, and travel speed can vary. The recent efforts \cite{bourque2019,berger2021near} and, largely, \cite{Royset2010} deal with homogeneous searchers where all these characteristics are identical across the searchers. We permit the target to camouflage according to a random process. Thus, the target not only follows a random trajectory but its appearance along the trajectory is also random. It might become undetectable for some time periods and this adds variability to the searchers' effective sensor performance at any point in time. To the best of our knowledge, this feature has not been modeled earlier in the literature.

We start in Section 2 by formulating the search problem under consideration. Section 3 considers the most general conditional target path models and presents two linearizations, a preprocessing technique, an outer-approximation method based on lazy-constraints, and numerical results. Section 4 turns to the more special, Markovian target path models and develops a linearization and three cutting plane algorithms, with supporting numerical results. The paper ends with conclusions in Section 5.

\section{Problem Formulation}

In this section, we describe the search problem and propose a generic model formulation.

\subsection{Searchers and the Target}\label{SUB21}

We consider $L$ classes of searchers with each class $l\in\cL=\{1, \ldots, L\}$ containing $J_l$ identical searchers. The set of time periods is $\cT_0 = \{0\}\cup\cT$ with $\cT = \{1, \dots,T\}$. The search for the target may take place during time periods $t\in\cT$. During a time period $t\in\cT_0$, each searcher occupies a state $s\in \cS = \{1, \dots, S\}$ or is in transit between states. When occupying a state $s$, a searcher of class $l$ may select to move to any state adjacent to $s$ as defined by the forward star $\cF_l(s)\subset \cS$. We also let $\cR_l(s)\subset\cS$ denote the reverse star of state $s$, which represents the set of states from which a searcher of class $l$ can reach state $s$ without transiting through any intermediate state. A searcher of class $l$ requires $d_{l,s,s'}\geq 1$ time periods to move from state $s$ to state $s'\in\cF_l(s)$ and to carry out search in state $s'$ for one time period. We refer to $d_{l,s,s'}$ as the {\em travel time} even though it also includes the subsequent search time and typically would have $d_{l,s,s}=1$ when the searcher remains in state $s$.

We prefer the term ``state'' over ``cell'' despite the latter being more common in the literature; see, e.g., \cite{Royset2010,berger2021near}. ``State'' highlights the vast number of modeling possibilities beyond searching an area discretized into grid cells. For example, the search may take place inside an underground mine, inside a ship, in a building, or in an urban environment. In such situations, it becomes especially important to allow for varying travel times $d_{l,s,s'}$ that sometimes could be much greater than one time period.

We let $X_{l,s,s',t}$ denote the number of searchers of class $l$ that occupy state $s$ in time period $t\in\cT_0$ and that move to state $s'$ next, and let $X$ denote the vector with components $X_{l,s,s',t}$, $l\in\cL$, $s,s'\in\cS$, and $t\in\cT_0$. We refer to $X$ as a {\em search plan}.
In addition to the conditions imposed by the forward and reverse stars, a search plan is constrained in three ways:\\

\state Initial State. There is a special, {\em initial state} $s_+ \in \cS$ from which all searchers start at time period $0$. It can abstractly represent geographically distinct bases for the different classes of searchers as the travel time $d_{l,s_+,s}$ from $s_+$ to any other state $s$ may depend on $l$. (Further fidelity regarding starting states for the various searchers is easily implemented, but omitted here for notational simplicity.) The reverse star $\cR_l(s_+) = \{s_+\}$, indicating that a searcher cannot return to the initial state after it departs. However, since $\cF_l(s_+)$ may contain $s_+$, a searcher could remain in the initial state for a number~of~periods.\\

\state Deconfliction. We permit at most $n_{s,t}$ searchers to be in state $s$ at time period $t$. This constraint is motivated by safety concerns related to collisions, but could also be helpful in preventing search plans that overly concentrate on a few states. Our modeling framework easily accommodates a variety of other deconflication constraints as well, but we omit the details.\\

\state Endurance and Terminal State. For each class $l$, there is an {\em endurance level} $\tau_l$ which is the number of periods a searcher of that class can be absent from $s_+$ and $s_-\in \cS$, the latter being the {\em terminal state}. It has the forward star $\cF_l(s_-) = \{s_-\}$, which means that a searcher in the terminal state will remain there indefinitely. As we see in the below formulation, travel time from $s_+$ to the first state looked at and travel time from the last state to $s_-$ are not counted against $\tau_l$. For example, suppose that $\cS = \{1, \dots, 5\}$, $s_+ = 1$, $s_- = 5$, and consider the forward stars $\cF_1(1) = \{1,2\}$, $\cF_1(2) = \{2,3\}$, $\cF_1(s) = \{s-1,s,s+1\}$ for $s = 3, 4$, $\cF_1(5) = \{5\}$  and the reverse stars $\cR_1(1) = \{1\}$, $\cR_1(s) = \{s-1,s,s+1\}$ for $s = 2, 3$, $\cR_1(4) = \{3,4\}$, $\cR_1(5) = \{4,5\}$. If $T = 6$ and $\tau_1 = 3$, then a feasible plan for searcher 1 is to sequentially visit the states $1, 1, 2, 3, 4, 5$ because the searcher is outside of the initial and terminal states for no more than $\tau_1 = 3$ time periods.\\

We consider one target. During a time period $t\in \cT$, the target is in a state $s_t\in \cS\setminus\{s_+, s_-\}$ while operating in one of two modes: it might be camouflaged at that time as indicated by $c_t = 1$ or it might not be camouflaged specified by $c_t = 0$. We observe that the target is barred from the initial and terminal states of the searchers. A {\em target path} is the vector $\omega = (\omega_1, \dots, \omega_T)$ with $\omega_t = (s_t, c_t) \in (\cS\setminus\{s_+, s_-\})\times \{0,1\}$ specifying the state $s_t$ and mode $c_t$ for the target in time period $t$. The probability that the target follows path $\omega$ is $q(\omega)$. We denote by $\Omega \subset ((\cS\setminus\{s_+, s_-\}) \times \{0,1\})^T$ the set of all target paths with positive probability. Thus, $\sum_{\omega \in \Omega} q(\omega) = 1$. We assume that these target paths and probabilities are known. Since we adopt a stochastic model for target movement, it becomes immaterial whether the target {\em wants} to be detected or not. The target simply selects one target path according to the probabilities $q(\omega), \omega\in \Omega$ and follows it without any ``intelligent'' behavior.

While we only explicitly consider a single target, it is conceptually straightforward to extend the following formulations to multiple targets by adopting expected number of unique targets detected or related metrics as objective function. Since this only affects the objective function with the decision variables remaining the same, we conjecture that computing times will largely be unchanged compared to the single-target case. We omit a detailed discussion and refer to \cite{Royset2010} for ideas in this direction.

\subsection{Sensors} \label{SUB22}

We assume that each searcher is equipped with one imperfect sensor. Each time period $t\in\cT$ in which a searcher occupies a state, the searcher's sensor takes one {\em look} at its current state. When a searcher is in transit between states, the sensor is inactive. If a searcher of class $l$ occupies state $s$ in time period $t$ and $s'$ is the searcher's previous state, then the probability that the searcher's look at the state during time period $t$ detects the target, given it is in that state and is not camouflaged, is $g_{l,s',s,t}\in [0,1)$. We refer to this probability as the {\it glimpse-detection probability}. We assume that the searchers' looks can be viewed as statistically independent attempts at detecting the target. Hence, given a search plan $X$ and target path $\omega$, the probability that no searcher detects the target during $\cT$ becomes:
\begin{subequations}
\begin{align*}
&\prod_{l\in \cL}\prod_{s\in \cS}\prod_{t\in \cT} \prod_{\substack{s'\in\cR_l(s) \\ t-d_{l,s',s}\geq 0}}
(1-g_{l,s',s,t})^{\zeta_{s,t}(\omega) X_{l,s',s,t-d_{l,s',s}}}\\
& = \exp \Bigg(-
\sum_{l\in\cL}\sum_{s\in \cS} \sum_{t\in \cT} \sum_{\substack{s'\in\cR_l(s) \\ t-d_{l,s',s}\geq 0}}
-\ln (1-g_{l,s',s,t}) \zeta_{s,t}(\omega) X_{l,s',s,t-d_{l,s',s}} \Bigg)\label{eqn:nodet2},
\end{align*}
\end{subequations}
where $\zeta_{s,t}(\omega) = 1$ if $\omega=(\omega_1, \dots, \omega_T)$ has $\omega_t = (s,0)$, and $\zeta_{s,t}(\omega) = 0$ otherwise. For given $l,s',s,t$, there are four possible reasons why
\[
(1-g_{l,s',s,t})^{\zeta_{s,t}(\omega) X_{l,s',s,t-d_{l,s',s}}}
\]
would become 1 and thus causing this particular factor to not reducing the probability of non-detection: (i) the glimpse-detection probability $g_{l,s',s,t}$ could be 0 representing an ineffective sensor under these circumstances. For example, $t$ might represent nighttime or a time period with poor weather. (ii) No searchers of class $l$ are present in state $s$ at time period $t$, while previously in $s'$, i.e., $X_{l,s',s,t-d_{l,s',s}} = 0$. (iii) The target is not in state $s$ at time $t$, which causes $\zeta_{s,t}(\omega) = 0$. (iv) The target is in state $s$ at time $t$ but is camouflaged, i.e., $\omega_t = (s,1)$, which again causes $\zeta_{s,t}(\omega) = 0$.

We refer to the term $-\ln (1-g_{l,s',s,t})$ as the {\em detection rate} for a searcher of class $l$ in state $s$ at time $t$ when it previously occupied state $s'$. Generally, these detection rates can vary with $l,s',s,t$ but we assume that one can identify a positive number $\alpha$ and nonnegative {\em integers} $\beta_{l,s',s,t}$, $l\in \cL, s,s'\in \cS, t\in \cT$, such that
\begin{equation}\label{eqn:alphadef}
\alpha\beta_{l,s',s,t} = -\ln (1-g_{l,s',s,t})~ \mbox{ for all }~ l\in \cL, s\in \cS, t\in \cT, s'\in \cR_l(s) \,\mbox{ with }\, t - d_{l,s',s} \geq 0.
\end{equation}
This is a minor assumption as each number in a finite collection of rational numbers can be written as the product of a common scalar and an integer. We refer to $\alpha$ as the {\em base detection rate}, while $\beta_{l,s',s,t}$ is the {\em rate modification factor}. The motivation for the assumption stems from the linearization approaches below; see also \cite{Royset2010} which mentions this possibility while leaving out the details. The complexity of a problem instance turns out to be closely related to the size of the integers $\beta_{l,s',s,t}$. If the sensors are identical across classes, states, and time periods, then one can set all rate modification factors to 1. To take advantage of this particular structure in the formulation below, we leverage the auxiliary decision variable
\[
Z_{l,s,t} = \sum_{\substack{s'\in\cR_l(s) \\ t-d_{l,s',s}\geq 0}} \beta_{l,s',s,t} X_{l,s',s,t-d_{l,s',s}},
\]
which represents the {\em search effort} allocated to state $s$ at time period $t$ by class $l$.

\subsection{SP Model} \label{SUB23}

We next state an MINLP that models the search problem under consideration. It goes beyond the formulations in \cite{bourque2019,berger2021near} by considering different classes of searchers, varying travel times, deconflication constraints, and endurance limits. It is motivated by a model in \cite{Royset2010}, but extends it by accounting for a camouflaging target and limited search endurance.
Table \ref{T5} provides a summary of the notation used.

\begin{table}[H]
	\centering
	\caption{Notation for model {\bf SP}}
	\label{T5}
	\setlength\arrayrulewidth{1.25pt}
		\begin{tabular}{ll}
		\hline
        \textbf{Indices} \\
        \hline
        $s,s',s_t$       & State: $s,s',s_t\in\cS=\{1,\ldots,S\}$\\
        $t,t'$           & Time period: $t,t'\in\cT_0=\{0\}\cup\cT$, $\cT = \{1,\ldots,T\}$ \\
        $l$              & Searcher class: $l\in\cL=\{1,\ldots,L\}$\\
        $c,c',c_t$          & Mode: $c = 1$ means camouflage; $c=0$ means no camouflage\\
        $\omega$         &  \multicolumn{1}{p{11cm}}{Target path: $\omega = (\omega_1, \dots, \omega_T)\in \Omega$, with $\omega_t = (s_t,c_t)\in (\cS\setminus\{s_+,s_-\}) \times \{0,1\}$}\\
        \hline
        \textbf{Sets} & \\
        \hline
        $\cF_l(s)\subseteq\cS$     & Forward star of state $s$ for searchers of class $l$ \\
        $\cR_l(s)\subseteq\cS$     & Reverse star of state $s$ for searchers of class $l$ \\
        \hline
        \textbf{Parameters}  \\
        \hline
        $\alpha$ & \multicolumn{1}{p{11cm}}{Base detection rate; positive real number}\\
        $\beta_{l,s',s,t}$ & \multicolumn{1}{p{11cm}}{Rate modification factor for a searcher of class $l$ while it occupies state $s$ in time period $t$ and $s'$ is its previous state; nonnegative integer}\\
        $\zeta_{s,t}(\omega)$ & \multicolumn{1}{p{11cm}}{1 if $\omega=(\omega_1, \dots, \omega_T)$ has $\omega_t = (s,0)$; zero otherwise}\\
        $s_{+} \in \cS$            & Initial state; $\cR_l(s_{+}) = \{s_{+}\}$\\
        $s_{-} \in \cS$            & Terminal state; $\cF_l(s_{-}) = \{s_{-}\}$\\
        $J_l$    & Number of searchers of class $l$; positive integer\\
        $q(\omega)$  & Probability of target path $\omega$; positive value with $\sum_{\omega \in \Omega} q(\omega) = 1$\\
        $d_{l,s,s'}$     & \multicolumn{1}{p{11cm}}{Number of time periods needed for a searcher of class $l$ to move directly from state $s$ to state $s'$ and search in $s'$; positive integer} \\
        $n_{s,t}$        & \multicolumn{1}{p{11cm}}{Maximum number of searchers in state $s$ at time period $t$; nonnegative integer} \\
        $\tau_l$         & Endurance of searchers of class $l$; positive integer\\
        $m_{l,s,t}$      & \multicolumn{1}{p{11cm}}{Maximum search effort from class $l$ in state $s$ at time period $t$; $m_{l,s,t} = \sum_{s'\in \cR_l(s):t-d_{l,s',s}} \beta_{l,s',s,t} \min\{J_l,n_{s,t}\}$} \\
        \hline
        \textbf{Decision Variables}  & \\
        \hline
        $X_{l,s,s',t}$ & \multicolumn{1}{p{11cm}}{Number of searchers of class $l$ in state $s$ at time period $t$ and that move to state $s'$ next; $X$ denotes the vector with components $X_{l,s,s',t}$, $l\in\cL, s,s'\in \cS, t\in\cT_0$}\\
        $Z_{l,s,t}$  &  \multicolumn{1}{p{11cm}}{Search effort from class $l$ in $s$ at time period $t$, $l\in \cL, s\in \cS, t\in \cT$; $Z$ denotes the vector with components $Z_{l,s,t}, l\in \cL, s\in \cS, t\in \cT$}\\
        $M_{l,t}$  &  \multicolumn{1}{p{11cm}}{Number of searchers of class $l$ that start their mission at time period $t$; $M$ denotes the vector with components $M_{l,t}$, $l\in \cL, t\in \cT$}\\
\end{tabular}
\end{table}

The MINLP takes the following form:
\begin{subequations}
\begin{align}
\textbf{SP:} \; \; \nnmin_{X,Z,M} & ~~f(Z) =\sum_{\omega\in\Omega} q(\omega)
   \exp\Bigg(-\sum_{l\in\cL} \sum_{\substack{s\in \cS \\ s\not\in\{s_+,s_-\}}} \sum_{t\in \cT} \zeta_{s,t}(\omega)\alpha Z_{l,s,t}\Bigg)
     \label{eqn:SPXobj}\\
\text{subject to }  & \sum_{\substack{s'\in\cR_l(s) \\ t-d_{l,s',s}\geq 0}} X_{l,s',s,t-d_{l,s',s}}  =  \sum_{s'\in\cF_l(s)}X_{l,s,s',t}, \; \;  l\in \cL, s\in \cS, t\in\cT\label{eqn:SPXflow}\\&  \sum_{s\in\cF_l(s_+)} X_{l,s_+,s,0}  =  J_{l}, \; \;  l\in \cL\label{eqn:SPXini}\\
&  \sum_{\substack{s\in \cF_l(s_+) \\ s\not\in \{s_{+},s_{-}\}}} X_{l,s_+,s,t} = M_{l,t}, \; \;  l\in \cL, t\in \cT_0\label{eqn:takeoff1}\\
&  \sum_{\substack{s \in \cS \\ s\not\in \{s_{+},s_{-}\}}}\sum_{s'\in\cF_l(s)} X_{l,s,s',t} \leq \sum_{t - \tau_l + 1 \leq t' \leq t} M_{l,t'}, \; \;  l\in \cL, t\in \cT_0\label{eqn:takeoff2} \\
&  \sum_{\substack{s'\in\cR_l(s) \\ t-d_{l,s',s}\geq 0}} \beta_{l,s',s,t} X_{l,s',s,t-d_{l,s',s}} = Z_{l,s,t}, \; \; l\in \cL, t\in \cT, s\in \cS\; \; \label{NEW1}\\
& \sum_{l\in \cL} \sum_{\substack{s'\in\cR_l(s) \\ t-d_{l,s',s}\geq 0}} X_{l,s',s,t-d_{l,s',s}} \leq n_{s,t}, \; \; t\in \cT, s\in \cS \label{NEW2}   \\
& X_{l,s,s',t} \in \big\{0, 1, 2, \dots, \min\{J_l,n_{s,t}\}\big\}, \; \;  l\in \cL,s,s'\in \cS,t \in \cT_0\label{REL-INTE}\\
& M_{l,t} \in \big\{0, 1, 2, \dots, \min\{J_l,n_{s_+,t}\}\big\}, \; \;  l\in \cL,t \in \cT_0\label{REL-INTEb}\\
& Z_{l,s,t} \in \{0, 1, 2, \dots, m_{l,s,t}\}, \; \; l \in \cL, t\in \cT, s\in \cS. \label{NEW3}
\end{align}
\end{subequations}

The objective function \eqref{eqn:SPXobj}, denoted by $f(Z)$, gives the probability of {\em not} detecting the target during $\cT$ and is obtained from the derivations in Subsection \ref{SUB22} by applying the total probability theorem. It leverages the auxiliary decision vector $Z$ assigned in \eqref{NEW1}. In view of \eqref{eqn:alphadef}, $\exp(\alpha Z_{l,s,t})$ gives the probability that class $l$ fails to detect the target  in state $s$ at time period $t$, given the target is there and it is not camouflaging.

Constraints (\ref{eqn:SPXflow}) and (\ref{eqn:SPXini}) enforce route continuity and define initial conditions for the searchers, respectively. The constraints \eqref{eqn:takeoff1} ensure that $M_{l,t}$ represents the number of searchers of class $l$ that moves away from the initial state in time period $t$, i.e., start their mission. The constraints \eqref{eqn:takeoff2} prevent searchers from being outside the initial and terminal states for more than $\tau_l$ time periods. Specifically, the right-hand side of \eqref{eqn:takeoff2} sums up the number of searchers of class $l$ that has started their mission during time periods $t, t-1, \dots, t-\tau_l + 1$. This number cannot be exceeded by the left-hand side of \eqref{eqn:takeoff2}, which gives the number of searchers of class $l$ on mission at time period $t$. Thus, searchers of class $l$ that started their mission prior to $t-\tau_l + 1$ cannot be in any other state than $s_-$. To the best of our knowledge, endurance constraints of this kind have not been considered earlier in the search theory literature.
Deconfliction constraints \eqref{NEW2} limit the number of searchers that can occupy a state in any time period. It can be adjusted in various ways such as being implemented for each class $l$ individually.

We can reduce the size of {\bf SP} by defining $Z_{s,t} = \sum_{l\in \cL} Z_{l,s,t}$, but the present formulation affords some simplifications. If each $\beta_{l,s',s,t} = 1$, then every $X_{l,s,s',t}$ can be relaxed to a continuous variable. This is not the case in a formulation with the aggregated variables $Z_{s,t}$.

{\bf SP} is a convex MINLP because its continuous relaxation has a convex nonlinear objective function and a polyhdedral feasible set. The difficulty of solving {\bf SP} depends on various parameters as examined below. The movement of the target between states and the switch in and out of camouflaging mode enter {\bf SP} only through the set of target paths $\Omega$, which are weighted according to the probabilities $q(\omega)$, $\omega\in \Omega$. Our formulation has the advantage that any (complicated) target path model can be considered, including non-Markovian models. It suffices to generate, ex-ante, the parameters $\zeta_{s,t}(\omega)$ for each path $\omega\in \Omega$. We refer to this most general setting as a {\em conditional target path model} and address~it~in~Section~\ref{SGM1}.

While conceptually simple, a conditional target path model might be computationally challenging to implement when the number of possible paths is large, i.e., the cardinality of $\Omega$ is large. A {\em Markovian target path model} affords a means to handle a massive number of target paths as we see in Section \ref{SGM2}.

\section{Conditional Target Paths}\label{SGM1}

In this section, we consider conditional target paths and thus make no assumptions about the stochastic model generating these paths beyond being able to compute ex-ante the parameters $\zeta_{s,t}(\omega)$. Subsection \ref{RF1} develops two equivalent linear models, a supporting preprocessing technique, and numerical results. Subsection \ref{AL1} presents an outer-approximation method based on lazy constraints, which improves computing times on difficult instances. Subsection \ref{sec:operational} discusses operational insights emerging from solving {\bf SP} in various settings.

\subsection{Linearization}\label{RF1}

The objective function \eqref{eqn:SPXobj} in {\bf SP} is a finite sum of the exponential function with arguments in the form of a sum of products of a nonnegative parameter by a bounded integer variable. It can therefore be linearized using additional variables and constraints \cite{Royset2010}. In addition to extending the linearization from \cite{Royset2010}, which deals with homogeneous searchers and no operational constraints, to the present setting, we also develop a novel linearization and a preprocessing technique.

The maximum search effort that the searchers collectively can muster across all time periods is
\[
N = \sum_{l\in \cL} \sum_{t\in \cT} \max_{s\in \cS\setminus\{s_+,s_-\}} m_{l,s,t} \ .
\]
Thus, the power in \eqref{eqn:SPXobj} cannot exceed $\alpha N$. A linearization of the exponential function needs to only cover the arguments $0$, $\alpha$, $2\alpha$, $\dots$, $\alpha N$.

We start by developing a new linearization by leveraging the fact that minimizing $\exp(-\alpha Y)$ over $Y\in \{0, 1, 2, \dots, N\} \cap \cY$, where $\cY$ represents constraints, is equivalent to the problem
\begin{equation}\label{eqn:upperapprox}
\nnmin_{Y \in \cY,W_0, \dots, W_N}  \sum_{i=0}^N W_i e^{-i\alpha} \text{ subject to } \sum\limits_{i=1}^N i \, W_i = Y, ~\sum\limits_{i=0}^N W_i = 1,~ W_i \in [0,1], ~i =0,1,2, \ldots,N.
\end{equation}
At optimality, each $W_i$ must take value 0 or 1 because the exponential function is strictly convex, which means that one can restrict $W_i$ to be binary from the outset. Replicating the process for each $\omega\in \Omega$ in the context of {\bf SP}, we reformulate {\bf SP} as the following mixed-integer linear program (MILP):
\begin{align}
\textbf{CSP-U:} \; \nnmin_{X,Z,W,M} & \; \sum_{\omega\in\Omega} q(\omega) \sum_{i=0}^N W_i(\omega) e^{-i\alpha} \notag\\
\text{subject to  } 	&\eqref{eqn:SPXflow}\mbox{-}\eqref{NEW3} \notag\\
& \sum\limits_{i=1}^N i \, W_i(\omega) = 	\sum_{l \in \cL} \sum_{\substack{s \in \cS \\ s\not\in \{s_{+},s_{-}\}}} \sum_{t \in \cT} \zeta_{s,t}(\omega) Z_{l,s,t} \, , \quad \omega \in \Omega \notag\\
&\sum\limits_{i=0}^N W_i(\omega) = 1, \quad \omega \in \Omega \label{E3}\\
&W_i(\omega) \in [0,1], \quad \omega \in \Omega, ~i = 0,1, 2, \ldots,N. \label{E4}
\end{align}
Here, we denote by $W$ the vector with components $W_i(\omega)$, $\omega \in \Omega$, $i =\{0,1,\ldots,N\}$. The first letter in {\bf CSP-U} refers to the \underline{c}onditional target model, while the last letter hints to the \underline{u}pper approximation of the exponential function underpinning \eqref{eqn:upperapprox}. Note that there is no approximation in the present setting; {\bf CSP-U} is equivalent to {\bf SP}.

We also extend a linearization from \cite{Royset2010}, which gives the following MILP reformulation of {\bf SP}:
\begin{align}
\textbf{CSP-L:} \; \nnmin_{X,Y,Z,M} & \; \sum_{\omega\in\Omega} q(\omega) Y(\omega) \notag\\
\text{subject to } &  \eqref{eqn:SPXflow}\mbox{-}\eqref{NEW3} \notag\\
& \; e^{-i\alpha}  (1+i-ie^{-\alpha}) -e^{-i \alpha} (1-e^{-\alpha} ) \sum_{l \in \cL} \sum_{\substack{s \in \cS \\ s\not\in \{s_{+},s_{-}\}}} \sum_{t \in \cT} \zeta_{s,t}(\omega) Z_{l,s,t}  \leq Y(\omega) \notag \\
& \;\hspace{7cm} \omega \in \Omega,~ i =0,1,2, \ldots,N-1. \label{eqn:SP1Lnew}
\end{align}
The vector $Y$ consists of the free variables $Y(\omega), \omega\in \Omega$ introduced in the reformulation. As explained in \cite{Royset2010}, the constraints \eqref{eqn:SP1Lnew} represent $N$ secant cuts that are valid at integer points of the exponential function; this is replicated for each $\omega \in \Omega$. The last letter in the name {\bf CSP-L} recalls that each cut represents a \underline{l}ower approximation of the objective function in {\bf SP}. {\bf CSP-L} amounts to an improvement over the model {\bf SP1-L} in \cite{Royset2010} by considering multiple searcher classes, eliminating $|\Omega|$ unnecessary secant cuts (effectively replacing $N$ by $N-1$ in \eqref{eqn:SP1Lnew}), and accounting for endurance and deconfliction.

The linearizations {\bf CSP-U} and {\bf CSP-L} are both equivalent to {\bf SP}. The former adds $|\Omega|(N+1)$ variables and $(2N+4)|\Omega|$ constraints, while the latter adds only $|\Omega|$ variables and $|\Omega|N$ constraints. However, the added constraints in {\bf CSP-U} are relatively simple; either variable bounds or equality constraints. In contrast, all the new constraints in {\bf CSP-L} are more challenging inequality constraints. Regardless, the role of $N$ is central, with lower values affording significant savings in model size. The planning horizon $T$ and the number of searchers drive up $N$. The same holds for situations with varying detection rates, which produce rate modification factors $\beta_{l,s',s,t}$ larger than one.

As is the case for {\bf SP}, if each $\beta_{l,s',s,t} = 1$ in {\bf CSP-U} and {\bf CSP-L}, then every $X_{l,s,s',t}$ can be relaxed to a continuous variable. When possible, we take advantage of this fact. (Testing not reported here indicates significant reduction in computing time when using this relaxation. The alternative relaxation with $Z$ continuous and $X$ integer is significantly slower, which probably stems from the fact that $X$ is a much larger vector than $Z$.)\\

\state Computational Tests. We compare {\bf CSP-U} and {\bf CSP-L} in a preliminary computational study based on instances from \cite{Royset2010}. For reference, we also examine the standard solvers Baron, Bonmin, and Knitro \cite{KronqvistBernalLundellGrossmann.19}. There is a single class of searchers with unlimited endurance looking for a target that cannot go into camouflage mode. We also omit the deconfliction restrictions \eqref{NEW2}. This implies that the variable vector $M$ and the constraints \eqref{eqn:takeoff1} and \eqref{eqn:takeoff2} are superfluous. The state space is built as a square grid of cells, with an additional state $s_+$ representing the initial location of the searchers. (A terminal state $s_-$ is unnecessary when the searchers have unlimited endurance.) For example, a 9-by-9 grid of cells produces 81+1= 82 states. At any time period $t$, a searcher in state $s$, corresponding to a particular grid cell, can move to the cell above, below, right, or left to $s$ in the grid and this becomes its next state. We call these four states as well as $s$ itself the adjacent states of $s$. Diagonal moves are not allowed. On the boundary of the square grid of cells some of these options are eliminated as needed. The adjacent states define the forward star set $\cF_l(s)$. The reverse star of $s$ is defined analogously. The travel times $d_{l,s,s'}$ are always set to 1. The initial state $s_+$ has the three boundary cells in the upper-left corner as its forward star. The glimpse detection probabilities are invariant so that $\beta_{l,s',s,t}=1$ for all $l,s,s',t$, with $\alpha = -3\ln(0.4)/J_1$; here $J_1$ is the number of searchers of the first (and only) class. This calibration of $\alpha$ follows \cite{Royset2010} and allows for comparison as the number of searchers varies.

The target paths are generated ex-ante as follows. The number of cells along each edge of the square grid of cells is an odd number, so the center cell in the square grid is well defined. This center cell is the initial position of the target. From one time period to the next, the target can stay idle or move to any of the adjacent cells according to a transition matrix with probabilities defined as follows. The probability that the target remains in the same state is 0.6, with the probability of moving to any of the adjacent states is equal (i.e., usually 0.1 except if the target is on the boundary of the square grid of cells). We randomly generate $|\Omega|$ target paths according to these probabilities and set $q(\omega) = 1/|\Omega|$.

These model instances and those in the following are not constructed in response to a particular application, but rather designed to challenge the algorithms. Current and future applications might involve many searchers in the form of inexpensive drones or a few manned aircraft. The number of states can also vary greatly. The search for smugglers in the Eastern Tropical Pacific Ocean might involve thousands of states, two aircraft, 72 hourly time periods, and half-a-dozen targets \cite{Riley.23}. However, after preprocessing and decoupling the various targets we obtain a state space and planning horizon aligned with what is considered in this paper. 

All the models in this paper are coded in Python 3.7 and solved with Gurobi 9.1 on a Linux machine, with Intel Core i7-6700 CPU 3.40GHz processors and 64 GB installed physical memory.  For each instance, the relative optimality tolerance is 0.0001, and we use one thread only. If this tolerance is not achieved after 900 seconds, we report the optimality gap at 900 seconds in brackets in the tables below. The relative optimality gap is calculated as the ratio of the difference between the best integer solution and the best lower bound to the best lower bound.

\begin{table}[H]
	\centering
	\caption{{\small For $S = 82$ states, $|\Omega| = 1000$ target paths, and varying numbers of searchers and time periods: Solution time (sec.) to relative optimality gap of 0.0001 or, if not reached in 900 seconds, relative optimality gap in brackets after 900 seconds. Asterisk indicates that runtime is reduced to 17 seconds if $W_i(\omega)$ is restricted to binary in {\bf CSP-U}; $\infty$ indicates that no bound is available.}}
	{\small\begin{tabular}{|c|r|r|r|r|r|r|r|r|r|r|}
			\hline
			& \multicolumn{5}{c|}{$J_1=3$} & \multicolumn{5}{c|}{$J_1=15$} \\
			\hline
			$T$	&Baron	&Bonmin	&Knitro	&{\bf CSP-L}&{\bf CSP-U}	 &	Baron	&	Bonmin	&	Knitro	&	{\bf CSP-L} & {\bf CSP-U}\\
			\hline
			7	&113	&	9	&	17	&0.1&0.2&2	     &	9	  &	5	&0.9	&0.6\\
			8	&3	    &	14	&	23	&0.3&0.3&3	     &	15	  &	2	&1	&1 \\
			9	&48	    &	64	&	49	&2	&1	&10	     &	81	  &	12	&5	&3 \\
			10	&120	&285	&140	&5	&3	&23	     &	147	  &	8	&25	&*63 \\
			11	&[0.0153]&[0.0040]&200	&12	&6	&273	 &	461	  &	263	&436	&220 \\
			12	&[0.0482]&[0.0789]&451	&37	&7	&877	 &[0.4342]&	161	&82	&24 \\
			13	&[0.0367]&$[\infty]$	&$[\infty]$	&22	&10	&[0.0124]&[5.3512]&	284	&[0.0023]	&104 \\
			14	&[0.0577]&$[\infty]$	&$[\infty]$	&79	&18	&[0.0090]&[9.1903]&	797	&[0.0108]	&98 \\
			15	&[0.3043]&$[\infty]$	&$[\infty]$	&110&90	&$[\infty]$	 &$[\infty]$	  &$[\infty]$	&582   &279\\ 		
			\hline
	\end{tabular}}
	\label{table:LvsU-Tvarying}
\end{table}

Table \ref{table:LvsU-Tvarying} compares the Bonmin, Knitro, and Baron solvers with {\bf CSP-L} and {\bf CSP-U}. Direct solution of {\bf SP} using Bonmin, Knitro, and Baron appears less competitive: {\bf CSP-L} is faster than all the three solvers on 14 out of 18 instances; {\bf CSP-U} is faster than all the three solvers on 17 out of 18 instances and solves all of them within the 900-second time limit. Baron, Bonmin, and Knitro solve only 10, 9, and 14 out of 18 instances, respectively. Their failures often involve having found no feasible integer solution as indicated by $[\infty]$ in the table. A comparison between our linearizations shows that the new version {\bf CSP-U} tends to outperform {\bf CSP-L}, which in the present setting essentially coincides with a linearization proposed in \cite{Royset2010}. On 16 or 17 of the 18 instances, {\bf CSP-U} solves quicker than {\bf CSP-L}. The tolerance is reached in no more than 279 seconds with {\bf CSP-U}, while two instances cannot be solved in 900 seconds with {\bf CSP-L}. The advantage of {\bf CSP-U} over {\bf CSP-L} is more pronounced for instances with more searchers ($J_1 = 15$) compared to fewer searchers ($J_1 = 3$). We obtain similar results (not reported in detail) for instances with up to 32000 targets paths and 226 states in seconds. Interestingly, the solution time is {\em not} consistently increasing with the number of target paths and states.

In some cases a binary restriction on $W_i(\omega)$ in \eqref{E4} can be beneficial from a computational point of view. (Recall from the discussion after \eqref{eqn:upperapprox} that these variables indeed are binary at optimality.) For example, the instance with $T=10$ solves in 17 seconds with $W_i(\omega) \in \{0,1\}$ and in 63 seconds with $W_i(\omega) \in [0,1]$.


\begin{table}[H]
	\centering
	\caption{{\small For $S=82$ states, $|\Omega|=1000$ target paths, and varying time periods and numbers of searchers: Solution time (sec.) to relative optimality gap of 0.0001 or, if not reached in 900 seconds, relative optimality gap in brackets after 900 seconds.  Asterisk and dagger indicate that runtime is reduced to 17 seconds and 23 seconds, respectively, if $W_i(\omega)$ is restricted to binary in {\bf CSP-U}.}}
    {\small\begin{tabular}{|c|r|r|r|r|}
		\hline
		& \multicolumn{2}{c|}{$T=10$} & \multicolumn{2}{c|}{$T=15$} \\
		\hline
		$J_1$ & \textbf{CSP-L} & \textbf{CSP-U} & \textbf{CSP-L}  & \textbf{CSP-U}\\
		\hline
		3  & 5   & 3  & 110       & 89\\
		4  & 6   & 2  & 537       & 21\\
		5  & 9   & 4  & 33        & 42\\
		6  & 9   & 3  & 312       & 149\\
		8  & 15  & 4  & 124       & 152\\
		10 & 16  & 14 & 266       & 209\\
		15 & 26  & *63& 594       & 379\\
		20 & 34  & 39 & [0.0030]  & 503\\
		30 & 49  & 52 & [0.0021]  & 42\\
		50 & 51  & 15 & [0.0012]  & $\dagger$57\\
		\hline
	\end{tabular}}
	\label{table:LvsU-JTvarying}
\end{table}


The solution time appears to be an increasing function of the length of the planning horizon as seen in Table \ref{table:LvsU-JTvarying}, and this is also largely consistent with Table \ref{table:LvsU-Tvarying}. The effect of more searchers on the computing time is less clear. Instances with many searchers in Table \ref{table:LvsU-JTvarying} solve surprisingly quickly. The superiority of the new linearization {\bf CSP-U} becomes increasingly visible as the number of searchers and the length of the planning horizon increase.
For the largest instances with $J_1 \geq 30$ and $T=15$, Table \ref{table:LvsU-JTvarying} shows solution times for {\bf CSP-U} in tens of seconds while {\bf CSP-L} fails to produce the required optimality gap in 900 seconds. {\bf CSP-U} can also be solved with binary restrictions for $W_i(\omega)$, which is usually slower, but for 10 out of 52 instances in Tables \ref{table:LvsU-Tvarying}-\ref{table:LvsU-JTvarying} binary restrictions are slightly faster. The tables ignore such potential further improvements for {\bf CSP-U} unless the times become less than half in which case the instances are marked with asterisk and dagger in the tables.\\

\state Preprocessing. The linearizations of {\bf SP} involve a significant lifting of the decision space; it grows linearly in the number of target paths $|\Omega|$. The additional $|\Omega| N$ constraints in {\bf CSP-L} are also problematic. As a result, {\bf CSP-U} and {\bf CSP-L} can become prohibitively large for instances with many target paths, time periods, searchers, and/or varying rate modification factors. This motivates us to derive a preprocessing techniques to eliminate integer variables that can be proven to take value 0 at an optimal solution of {\bf CSP-U} or {\bf CSP-L} and to eliminate constraints that can be proven to be redundant.

If it can be determined a priori that no detection is possible in state $s$ during time period $t$, then {\it some of} the decision variables corresponding to the tuple $(s,t)$ can be fixed and/or removed. For this purpose we define the set $\mathcal D$ that includes all tuples $(s,t)$ for which detection is possible:
\vspace{-0.05in}
\begin{equation*}
	\label{SET}
	\mathcal D = \Big\{(s,t) \in \cS \times \cT ~\Big|~ \sum_{\omega \in \Omega} \zeta_{s,t}(\omega)>0\Big\}.
\end{equation*}
Let $\mathcal D^c$ denote the complement of $\mathcal D$. It follows that, if $(s,t) \in \mathcal D^c$,  having $Z_{l,s,t}>0$ will not reduce the probability of non-detection compared to having $Z_{l,s,t}=0$. Therefore, the corresponding integer variables $Z_{l,s,t}, (s,t) \in \mathcal D^c, l \in \cL$ can be removed from the formulation. Using this preprocessing approach, we obtain the following reduced-size formulations {\bf CSP-U-Pre} and {\bf CSP-L-Pre} for {\bf CSP-U} and {\bf CSP-L}, respectively:
\vspace{-0.05in}
\begin{subequations}
	\begin{align}
		\textbf{CSP-U-Pre:} \;
		\nnmin_{X,W,Z,M} & \; \sum_{\omega\in\Omega} q(\omega) \sum_{i=0}^N W_i(\omega) e^{-i\alpha} && \notag \\
		\text{subject to }  & \;
		\eqref{eqn:SPXflow}\mbox{-}\eqref{eqn:takeoff2}; \eqref{NEW2}\mbox{-}\eqref{REL-INTEb}; \eqref{E3}\mbox{-}\eqref{E4}  \notag \\
		& \; \sum\limits_{i=1}^N i \, W_i(\omega) = \sum_{l \in \cL} \sum_{(s,t) \in \mathcal D} \zeta_{s,t}(\omega) Z_{l,s,t}, \quad \omega \in \Omega \label{E2-NEW} \\
		&  \sum_{\substack{s'\in\cR_l(s) \\ t-d_{l,s',s}\geq 0}} \beta_{l,s',s,t} X_{l,s',s,t-d_{l,s',s}} = Z_{l,s,t}, \; \; l\in \cL, (s,t) \in \cD\; \; \label{NEW1-1}\\
		& Z_{l,s,t} \in \{0, 1, 2, \dots, m_{l,s,t}\}, \; \; l \in \cL, (s,t)\in \cD. \label{NEW3-1}
	\end{align}
\end{subequations}

\vspace{-0.05in}
\begin{subequations}
	\begin{align}
		\textbf{CSP-L-Pre:} \;
		\nnmin_{X,Y,Z,M} & \; \sum_{\omega\in \Omega} q(\omega) Y(\omega) &&\label{OBJ1}\\
		\text{subject to }  & \;
		\eqref{eqn:SPXflow}\mbox{-}\eqref{eqn:takeoff2}; \eqref{NEW2}\mbox{-}\eqref{REL-INTEb}; \eqref{NEW1-1}\mbox{-}\eqref{NEW3-1} \notag \\
		& \; e^{-i\alpha}  (1+i-ie^{-\alpha}) -e^{-i\alpha} (1-e^{-\alpha} ) \sum_{l \in \cL} \sum_{(s,t) \in \cD}  \zeta_{s,t}(\omega)  Z_{l,s,t} \leq Y(\omega) \notag \\
		& \;\hspace{7cm} \omega \in \Omega, i =\{0,1,\ldots,N-1\}. \label{eqn:SP1Llin-1}
	\end{align}
\end{subequations}

The preprocessing potentially reduces the size of the decision and constraint spaces in both {\bf CSP-U-Pre} and {\bf CSP-L-Pre}, and eliminates many vacuous constraints that otherwise would have entered \eqref{eqn:SP1Llin-1}. Numerical results comparing the efficiency of the formulations are provided next.

\subsection{Outer-Approximation Method} \label{AL1}

In this subsection, we develop an outer-approximation method {\tt OA} for solving large-scale instances of {\bf CSP-L-Pre} (and {\bf CSP-L}). An analogous approach for {\bf CSP-U} and {\bf CSP-U-Pre} is not possible.
While the preprocessing technique presented above provides a more compact reformulation, it remains nonetheless that the number of constraints \eqref{eqn:SP1Llin-1} can be extremely large. However, the vast majority of these constraints are not binding at an optimal solution.

The outer approximation outlined next builds on this observation and identifies a priori a vast set of constraints \eqref{eqn:SP1Llin-1} that are unlikely to impact the optimal solution, and can be viewed as \textit{lazy constraints} \cite{kleinert2021,lundell2019} and are defined as such in our algorithmic approach. They are at first removed from the formulation, giving a mixed-integer linear outer approximation (relaxation) {\bf OA$^{0}$} of problem {\bf CSP-L-Pre} (or {\bf CSP-L}) at the root node $0$ of the branch-and-bound (B\&B) tree.
Subsequently, at each node of the B\&B tree, we check whether the optimal solution at the current node  violates any such constraints. If so, the current optimal solution is discarded and the violated constraints are introduced in the updated outer approximation of all open nodes.
In short, the lazy constraints are moved to a pool and are initially removed from the constraint set before being (possibly) iteratively reinstated on an as-needed basis.
Caution must be exerted when selecting the lazy constraints and one should not be too aggressive. Indeed, the verification of whether a lazy constraint is violated is carried out each time a new incumbent solution is found and the overhead consecutive to the reinsertion of lazy constraints in the constraint set can be significant.

The challenge is to identify the constraints that can be removed so that (i) the size of the constraint set is reduced as much as possible and (ii) that few, if any, of the removed constraints will need to be reincorporated. For \eqref{eqn:SP1Llin-1}, we identify the levels of search effort that can be expected and this leads to an initial set of lazy constraints $\cL^0$:
\vspace{-0.05in}
\begin{equation}
	\label{LAZY}
	\mathcal{L}^0=  \left\{ \eqref{eqn:SP1Llin-1}: i \in \{0,\ldots, b_1 \} \cup \{b_2+1, b_2+2, \ldots, N\}  \ , \omega \in \Omega \right\} \ .
\end{equation}
The set $\mathcal L^0$ includes the constraints \eqref{eqn:SP1Llin-1} associated with an unlikely low and high number of looks as defined by the positive constants $b_1 < b_2 < N$.

We adopt the following notation. Let $\mathcal O$ denote the set of open nodes in the tree. Let $\mathcal{F}$ be the entire constraint set of problem {\bf CSP-L-Pre}, $\mathcal{L}^k$
be the set of lazy constraints at node $k$, $\mathcal{V}_L^k$ be the set of violated lazy constraints at $k$, and  $\mathcal{A}^k:= \cF \setminus \mathcal{L}^k$ be the set of active constraints at $k$, i.e., constraints included in the outer approximation considered at node $k$.

This leads to the outer-approximation method {\tt OA}: At the root node ($k=0$), we have $\mathcal{L}^0$  as defined in \eqref{LAZY}, $\mathcal{A}^0:= \mathcal F \setminus \mathcal{L}^0$, and $\mathcal{V}^0_L:= \emptyset$. At any node $k$, we solve the outer approximation
\vspace{-0.085in}
\[
\textbf{OA}^k: \; \nnmin \eqref{OBJ1} \quad \text{subject to} \quad (X,Y,Z,M) \in \mathcal{A}^k.
\vspace{-0.05in}
\]
Two cases exist for the optimal solution  $Z^{k*}$ of the continuous relaxation of $\textbf{OA}^k$:

\begin{enumerate}
	\item If  $Z^{k^*}$ is fractional, we introduce branching linear inequalities to cut off the fractional nodal optimal solution and continue the B\&B process.
	\item If $Z^{k^*}$ is integral and improves upon the incumbent solution, we check for possible violation of any lazy constraints.
	If any constraint in $\mathcal{L}^k$ is violated by    $Z^{k^*}$, we insert each  constraint violated in $\mathcal{V}_L^k \subseteq \mathcal{L}^k$ and discard $Z^{k*}$.
	We update the lazy and active constraint sets of each open node $o$ by letting $\mathcal{L}^o \leftarrow \mathcal{L}^o  \setminus \mathcal{V}_L^k$ and $\mathcal{A}^o \leftarrow \mathcal{A}^o  \cup \mathcal{V}_L^k$. On the other hand, if no lazy constraint in $\mathcal{L}^k$ is violated, $Z^{k*}$ becomes the incumbent solution and the node is pruned.
\end{enumerate}
In summary, the {\tt OA} method solves a reduced-size relaxation of {\bf CSP-L-Pre} at each node of the tree.
Each time {\bf OA}$^k$ provides an integral solution with better objective value than the incumbent solution, a verification is made if any lazy constraint is violated.  If it is the case, the incumbent integer solution is discarded and the violated lazy constraints are (re)introduced in the constraint set of all unprocessed nodes of the tree, thereby cutting off the current solution.  The above process terminates when all nodes are pruned.
We note that the callback verification is not performed at each node of the tree, but only when a better integer-valued feasible solution is found at a node.\\

\state Computational Tests. We next examine the efficiency of the {\tt OA} method and the effect of preprocessing as specified by {\bf CSP-U-Pre} and {\bf CSP-L-Pre}. Table \ref{table:LvsUvsOA-Jvarying} shows computing times for large instances generated as described in Subsection \ref{RF1}. (The table has occasional overlap with earlier tables and any discrepancy in the reported times are due to differences among randomly generated instances.) The preprocessing technique is typically beneficial, especially {\bf CSP-L-Pre} is an improvement over {\bf CSP-L}. {\bf CSP-U-Pre} is less consistent and might even add computing time compared to {\bf CSP-U} for instances when there are many searchers. In part, this is caused by the remarkable efficiency of {\bf CSP-U} on such instances. Generally, the best solution method appears to be the {\tt OA} method, which solves to optimality all instances in the allotted time and is the fastest on all but two instances. On the instances in Table \ref{table:LvsUvsOA-Jvarying}, {\bf CSP-L} is essentially identical to the approach proposed in \cite{Royset2010} but here falls behind with an order of magnitude longer computing times compared to the new approaches develop in the present paper.

\begin{table}[H]
	\centering
	\caption{{\small For $S=82$ states, $T = 15$ time periods, $|\Omega|=1000$ target paths, and varying numbers of searchers: Solution time (sec.) to relative optimality gap of 0.0001 or, if not reached in 900 seconds, relative optimality gap in brackets after 900 seconds.  Asterisk indicates that runtime is reduced to 23 seconds if $W_i(\omega)$ is restricted to binary in {\bf CSP-U}.}}
	{\small\begin{tabular}{|c|r|r|r|r|r|}
		\hline
		$J_1$ & \textbf{CSP-L} & \textbf{CSP-L-Pre} & \textbf{CSP-U} & \textbf{CSP-U-Pre} & {\tt OA} Method\\
		\hline
		3  & 110      & 60       & 89  & 89      &  54\\
		4  & 537      & 20       & 21  & 30      &  11\\
		5  & 33       & 20       & 42  & 27      &  9\\
		6  & 312      & 95       & 149 & 83      &  39\\
		8  & 124      & 28       & 152 & 30      &  11\\
		10 & 266      & 83       & 209 & 108     &  28\\
		15 & 594      & 313      & 379 & [0.0002]&  112\\
		20 & [0.0030] & 553      & 503 & 421     &  156\\
		30 & [0.0021] & [0.0002] & 42  & 206     & 271\\
		50 & [0.0012] & [0.0001] & *57 & 695     & 831\\
		\hline
	\end{tabular}}
	\label{table:LvsUvsOA-Jvarying}
\end{table}

Next, we consider more complex instances with a camouflaging target and searchers from two classes varying in their endurance level, which then activates constraints \eqref{eqn:takeoff1} and \eqref{eqn:takeoff2}. (We still omit deconflication constraints \eqref{NEW2}, which can be operationally important but produce simpler instances as many suboptimal search plans are immediately ruled out.)  The states are generated from a square grid of cells with an additional initial state as earlier, but now there is also a terminal state $s_-$. The reverse star of $s_-$ consists all the states corresponding to the bottom row of cells in the grid. The searchers otherwise move as earlier. The endurance of the searchers in class 1 and 2 is $\lfloor 0.8T \rfloor$ and $\lfloor 0.6T \rfloor$, respectively. For a total number of searchers $J$, the number of searchers in class 2 is $J_2=\lfloor 0.7 J \rfloor$, while the number of searchers in  class 1 is $J_1 = J - J_2$.
From a current state $s$, the target can opt to move to an adjacent state as before or to stay idle and transition into camouflage model. Once the target enters camouflage mode, it must stay in the same state in the next period, either camouflaged or not.
The target transition probabilities between states are as follows.
If occupying a state $s$ in non-camouflage mode, the target moves into camouflage mode (in the same state) with probability $0.1$. Otherwise, the target stays in $s$ in non-camouflage mode with probability $0.5$ or moves to an adjacent state with equal probability. When camouflaged, regardless of state, the target remains in camouflage mode with probability $1/6$ and comes out of it with probability $5/6$. Following these probabilities, we generate ex-ante $|\Omega|$ target paths.

Table \ref{table:CamoEndurance} summarizes the computing times for these instances across the various approaches. {\bf CSP-U} retains an edge over {\bf CSP-L} for instances involving less that 10 searchers. Interesting, the preprocessing technique delivers inconsistently on these instances, possibly due to the added complexity caused by the endurance constraints.  The best solution method appears to be the {\tt OA} method, which solves to optimality all instances in the allotted time and is always the fastest.
The solution time with the {\tt OA} method is not adversely affected by an increase in the number of searchers. The instance with 50 searchers can, for example, be solved about 40\% quicker than the one comprising 8 searchers.

\begin{table}[H]
	\centering
	\caption{{\small For $S=83$ states, $T = 15$ time periods, $|\Omega|=1000$ target paths, camouflaging target, and varying numbers of searchers $J$ split between two classes: Solution time (sec.) to relative optimality gap of 0.0001 or, if not reached in 900 seconds, relative optimality gap in brackets after 900 seconds.}}
	{\small\begin{tabular}{|c|r|r|r|r|r|}
		\hline
		$J$ & \textbf{CSP-L} & \textbf{CSP-L-Pre} & \textbf{CSP-U} & \textbf{CSP-U-Pre} & {\tt OA} Method\\
		\hline
		3  & 62  & 125 & 45  & 46   & 10 \\
		4  & 239 & 169 & 129 & 167  & 19 \\
		5  & 105 & 224 & 72  & 63   & 19 \\
		6  & 156 & 232 & 59  & 111  & 46 \\
		8  & 258 & 397 & 308 & 318  & 105 \\
		10 & 56  & 39  & 76  & 76   & 36 \\
		15 & 209 & 203 & 275 & 183  & 43 \\
		20 & 177 & 142 & 230 & 369  & 103 \\
		30 & 590 & 398 & [0.0002]   & [0.0002] & 125 \\
		50 & 155 & 173 & 234 & 144  & 64 \\
		\hline
	\end{tabular}}
	\label{table:CamoEndurance}
\end{table}

\subsection{Operational Insights}\label{sec:operational}

{\bf SP} enables an analyst or autonomous system to consider many different factors during the planning of a search mission. Next, we discuss the operational impact of limited endurance and varying travel times. We also quantify the difference between having many poor searchers compared to a few good ones.\\

\state Endurance and Travel Time. In an instance with 5 searchers, 2 in class 1 with endurance 12 and 3 in class 2 with endurance 9, we consider a 9-by-9 grid producing $S = 83$ states including the initial and terminal states as earlier. The planning horizon is $T =15$. We construct $|\Omega|=1000$ target paths without using the camouflage options as described in Subsection \ref{RF1}. The detection rate is the same for both classes, so $\beta_{l,s',s',t} = 1$ and $\alpha = -3\ln(0.4)/J$, where $J = 2 + 3 = 5$. Table \ref{table:path} shows an optimal search plan with objective function value 0.4244 using row-column notation to specify the state of each searcher during a time period. For example, the first searcher in class 1 stays in the initial state $s_+$ for three periods before moving to row 4, column 1 in the 9-by-9 grid as indicated by the pair $(4,1)$ in Table \ref{table:path}. In fact, $s_+$ has only this state in its forward star as hinted to in the table where every searcher moves to $(4,1)$ after departing $s_+$. We recall that the target starts in the middle of the grid: row 5, column 5. The search plan is thus meaningful with the searchers starting on the left rim and moving right as time progresses. In the absence of endurance constraints, all the searchers would obviously prefer to initialize their mission immediately. However, Table \ref{table:path} shows the interesting effect that under endurance limitations it is better for most of the searchers to wait a number of time periods and let the target ``come to them.'' The first time a searcher can encounter the target is in state $(4,3)$ in time period 3. But in periods 4 and 5, the target might have reached as far west as columns 2 and 1, respectively. Thus, a searcher starting its mission in period 5 or later may detect the target on its first look. The endurance constraints \eqref{eqn:takeoff1} and \eqref{eqn:takeoff2} introduce a delicate trade-off between searching early while the target is ``concentrated'' in the center of the grid cells but facing more ``wasted'' travel time versus searching late with the target being closer but more dispersed. For the present instance, the reverse star of $s_-$ consists all the states corresponding to the bottom row of cells in the grid, which we see the second searcher from class 1 moves toward as the time progresses. The other searchers remain on mission as we reach the planning horizon thus avoid having to enter $s_-$. This end-of-planning-horizon effect can be adjusted as needed with slight modification of constraints in {\bf SP}.

\begin{table}[H]
	\centering
	\caption{{\small Optimal search plan for $S=83$ states, $T = 15$ time periods, $|\Omega|=1000$ target paths, and 5 endurance-constrained searchers.}}
	\begin{small}
		\begin{tabular}{|c|c|c|c|c|c|c|c|c|c|c|c|c|c|c|c|}
			\hline
			& \multicolumn{15}{c|}{Time period $t$}\\
			\hline
			Class $l$		& 1       & 2       & 3       & 4       & 5       & 6       & 7     & 8     & 9     & 10    & 11    & 12    & 13     & 14     & 15     \\ \hline
			1 & $s_+$ & $s_+$& $s_+$ &  4,1    &  4,2    &  4,3    &  4,4    &  4,5    &  5,5    &  5,6  &  5,5  & 5,6  &  6,6  &  6,7  &  6,6    \\ \hline
			1 &  4,1    &  4,2    &  4,3    &  4,4    &  5,4    &  5,5    &  5,5  &  6,5  &  7,5  &  8,5  &  9,5  &  9,4  & $s_-$ & $s_-$ & $s_-$ \\ \hline
			2 & $s_+$ & $s_+$& $s_+$ & $s_+$ & $s_+$ & $s_+$ &  4,1  &  4,2  &  4,3  &  4,4  &  4,5  &  4,6  & 5,6   &  5,7   &  5,8   \\ \hline
			2 & $s_+$ & $s_+$& $s_+$ & $s_+$ & $s_+$ & $s_+$ &  4,1  &  4,2  &  4,3  &  4,4  &  4,5  &  5,5  & 5,5   &  5,5   &  5,6   \\ \hline
			2 & $s_+$ & $s_+$& $s_+$ & $s_+$ & $s_+$ & $s_+$ &  4,1  &  4,2  &  4,3  &  5,3  &  5,4  &  5,4  & 4,4   &  4,5   &  4,6   \\ \hline
		\end{tabular}
	\end{small}
	\label{table:path}
\end{table}

To illustrate the effect of other forward/reverse stars and travel times, which up to now has consisted of one-cell steps with $d_{l,s,s'}=1$, we slightly modify the instance by splitting class 1 into two classes: 1A and 1B, each with one searcher. The searcher in class 1A has augmented forward and reverse stars. In addition to the five states (stay, one cell up, one cell down, one cell left, one cell right) presently considered, we add the four states two cells up, two cells down, two cells left, and two cells right, again omitting nonexisting states outside the 9-by-9 grid of cells. Regardless, the travel time is $d_{l,s,s'}=1$. This means that the searcher is (potentially) faster than the searcher of class 1B, which retains the earlier forward/reverse star.
We split class 2 into three classes: 2A, 2B, and 2C, each with one searcher. The searcher in class 2A has augmented forward and reverse stars as class 1A.  The searcher in class 1B has the augmented forward and reverse stars as 1A, but the travel time is $d_{l,s,s'}=2$ if the searcher moves two cells, and otherwise $d_{l,s,s'}=1$. The searcher in class 2C is regular as for class 1B. Table \ref{table:path2} shows an optimal search plan with objective function value 0.4067. The improvement in probability of detection as compared to the search plan in Table \ref{table:path} stems from the faster searchers of class 1A and 2A; they move quickly toward the center of the grid cells. The searcher of class 2B has additional flexibility compared to Table \ref{table:path}, but does not leverage it because moving two cells in two periods without a look in the first cell cannot be better than moving one cell in one time period and then moving another cell in another time period while looking in both.

\begin{table}[H]
	\centering
	\caption{{\small Optimal search plan for $S=83$ states, $T = 15$ time periods, $|\Omega|=1000$ target paths, and 5 endurance-constrained searchers with varying forward/reverse stars and travel times.}}
	\begin{small}
		\begin{tabular}{|c|c|c|c|c|c|c|c|c|c|c|c|c|c|c|c|}
			\hline
			& \multicolumn{15}{c|}{Time period $t$}\\
			\hline
			Class $l$		& 1       & 2       & 3       & 4       & 5       & 6       & 7     & 8     & 9     & 10    & 11    & 12    & 13     & 14     & 15     \\ \hline
			1A & $s_+$ & $s_+$& $s_+$ & 4,1   & 4,3   & 4,5   & 5,5   & 5,4   & 5,6   & 5,5 & 6,5 & 6,4 & 6,6 & 6,5 & 6,6   \\ \hline
			1B & 4,1   & 4,2   & 4,3   & 4,4   & 5,4   & 5,5   & 5,5 & 6,5 & 7,5 & 8,5 & 9,5 & 9,4 & $s_-$ & $s_-$ & $s_-$ \\ \hline
			2A & $s_+$ & $s_+$& $s_+$ & $s_+$ & $s_+$ & $s_+$ & 4,1 & 4,3 & 4,5 & 5,5 & 5,6 & 5,4 & 5,5  & 5,7  & 5,8  \\ \hline
			2B & $s_+$ & $s_+$& $s_+$ & $s_+$ & $s_+$ & $s_+$ & 4,1 & 4,2 & 4,3 & 4,4 & 5,4 & 5,5 & 5,5  & 5,6  & 5,6  \\ \hline
			2C & $s_+$ & $s_+$& $s_+$ & $s_+$ & $s_+$ & $s_+$ & 4,1 & 4,2 & 4,3 & 5,3 & 5,4 & 5,5 & 4,5  & 4,5  & 4,6  \\ \hline
		\end{tabular}
	\end{small}
	\label{table:path2}
\end{table}

\vspace{0.05in}

\state Camouflage and Sensor Quality. We return to the setting at the end of Subsection \ref{AL1} and Table \ref{table:CamoEndurance}: there are two classes of searchers subject to endurance constraints and a target with camouflaging capability. As before, we use $\beta_{l,s',s',t} = 1$ and $\alpha = -3\ln(0.4)/J$, where $J = J_1 + J_2$. We examine the choice between acquiring many inexpensive but poor searchers or adopting few effective searchers at a higher cost. Our model of $\alpha$ as a function of the number of searchers $J$ has the consequence that $\alpha J$ is a constant. Thus, the power that can be mustered in the objective \eqref{eqn:SPXobj} is the same regardless of $J$. This means that having 10 searchers is in this sense equivalent to have 20 searchers because the former has an $\alpha$ twice as large as that of the latter. If each of the 10 more capable searchers are twice as expensive as each of the 20 less capable ones, then one might be indifferent between choosing 10 good versus choosing 20 poor searchers. The middle two rows, second column, in Table \ref{table:Minvalue} show that the objective value for the optimal search plans in these cases are indeed close:  0.4639 versus 0.4613. However, the slight detection improvement in the case of 20 searchers is not a coincidence. The case with 20 poor searchers produces a relaxation of {\bf SP} compared to the case with 10 good searchers because, in the absence of the deconflication constraint \eqref{NEW2}, the 20 poor searchers can always pair up to make a ``double-searcher'' of the same quality as any of the 10 good searchers. Going from 10 to 20 searchers, the change is minor but becomes more prevalent when we consider 50 searchers; see last row of Table \ref{table:Minvalue}. The effect appears to be reversed when we compare 5 and 10 searchers. However, the 10-searcher case is not a relaxation of the 5-searcher case because the latter has 2 searchers with 12-time-period endurance searchers and 3 searchers with 9-time-period endurance, while the former has 3 and 7 searchers for the two classes. Thus, the 10-searcher case has a slight endurance disadvantage and this causes the objective function value to increase. We also report the computing times for two methods in columns 3-4 of Table \ref{table:Minvalue}.

\begin{figure} [h!]
	\begin{center} \underline{Period 15:  Optimal searcher location at period $t=15$ with $J=50$ searchers}
	\end{center}
	\begin{minipage}{.48\linewidth}
		\vspace{-0.1in}
		\includegraphics[width=5.6cm, height = 5.6cm]{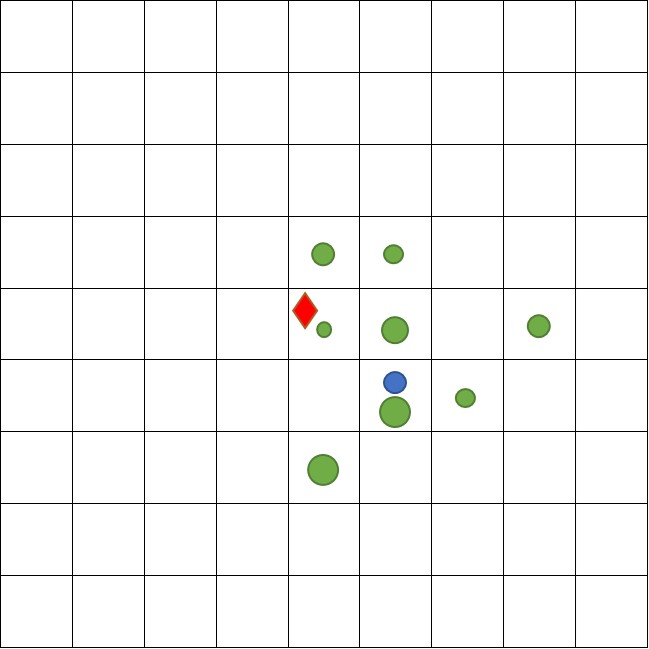}
		\centering
		\caption{{\bf With camouflage.} For class 1 (blue): 4 and 11 searcher in state (row, column) (6,6) and $s_-$, respectively. For class 2 (green): 4, 3, 2, 5, 4, 6, 3, 6, and 2 searchers in state
    (4,5), (4,6), (5,5), (5,6), (5,8), (6,6), (6,7), (7,5), and $s_-$.
    		}
		\label{fig-cam-50}
	\end{minipage}
	\hfill \hfill
	\begin{minipage}{.48\linewidth}
		\vspace{-0.1in}
		\includegraphics[width=5.6cm, height = 5.6cm]{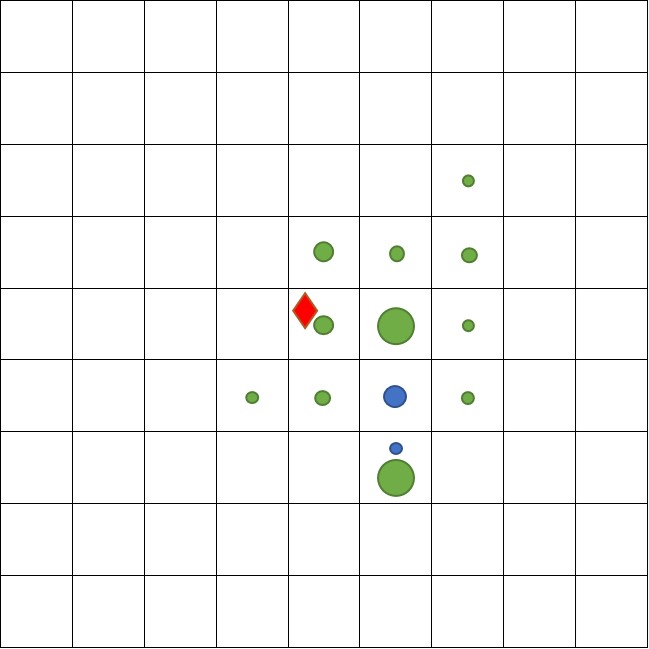}
		\centering
		\caption{{\bf Without camouflage.} For class 1 (blue): 4, 1, and 10 searchers in state (row, column) (6,6), (7,6), and $s_-$, respectively. For class 2 (green): 1, 3, 2, 2, 3, 7, 1, 1, 2, 1, 7, and 5 searchers in state
    (3,7), (4,5), (4,6), (4,7), (5,5), (5,6), (5,7), (6,4), (6,5), (6,7), (7,6), and $s_-$.
		}
		\label{fig-noncam-50}
	\end{minipage}%
\end{figure}

\begin{table}[H]
	\centering
	\caption{{\small For $S=83$ states, $T = 15$ time periods, $|\Omega|=1000$ target paths, and varying numbers of searchers and camouflaging capability: Min-value and solution time (sec.) to relative optimality gap of 0.0001 or, if not reached in 900 seconds, relative optimality gap in brackets after 900 seconds.}}
	{\small\begin{tabular}{|c|r|r|r|r|r|r|}
		\hline
		& \multicolumn{3}{c}{Camouflage} & \multicolumn{3}{|c|}{No camouflage}\\
		\hline
		$J$  & min-value & \textbf{CSP-L-Pre} & {\tt OA} method  & min-value & \textbf{CSP-L-Pre} & {\tt OA method}\\
		\hline
		5    & 0.4561 & 77  & 42   & 0.3500& [0.0016] & 321\\
		10   & 0.4639 & 78  & 31   & 0.3419& 229      & 76\\
		20   & 0.4613 & 509 & 315  & 0.3404& [0.0003] & 513\\
		50   & 0.4600 & 482 & 244  & 0.3399& [0.0001] & 134\\
		\hline
	\end{tabular}}
	\label{table:Minvalue}
\end{table}
\vspace{-0.085in}
We repeat the above calculations for a target that moves without camouflaging as described in Subsection \ref{RF1}; see the last three columns of Table \ref{table:Minvalue}. The probability of detecting the target improves with 0.10-0.13 because now the target can be detected everywhere along its path. We observe that the computing times for both the {\tt OA} method and {\bf CSP-L-Pre} tends to be less when the target can use camouflage. This is caused by a tighter concentration of likely target locations in the case of camouflage; it becomes less mobile with our parameter settings and the searchers' have fewer meaningful choices.
Figures \ref{fig-cam-50} and \ref{fig-noncam-50} illustrate the location of the 50 searchers from the last row of Table \ref{table:Minvalue} at time period 15. Here, the radius of a circle is proportional to the number of searchers occupying the corresponding state. The diamond indicates initial location for the moving target. Blue and green circles represent class 1 and class 2, respectively.
Figure \ref{fig-noncam-50} shows a wider spread of the searchers in the absence of camouflaging as compared to searchers concentrating on a less mobile, camouflaging target in Figure \ref{fig-cam-50}. At time period 15, the searchers tend to be on the eastern side as they have ``cleared'' the western side after entering at row 4, column 1.

\section{Markovian Target Paths}\label{SGM2}

We next present results for {\bf SP} under the assumption that the target moves according to a Markov chain, which thus defines the target paths $\Omega$ and the associated probabilities $q(\omega)$ by Markov transition matrices. Subsection \ref{RF2} presents a linear reformulation and Subsection \ref{AL2} develops three cutting plane methods. Numerical results appear in Subsection \ref{NT2}.

\subsection{Linearization} \label{RF2}

While the linearizations {\bf CSP-U} and {\bf CSP-L} remain valid for Markovian target paths, they tend to become prohibitively large unless the underlying state transition matrices are sparse or one adopts a sample average approximation with few sampled target paths. As noted by \cite{Royset2010} and refined in \cite{berger2021near}, the Markov structure affords an alternative linearization approach. These earlier studies focus on homogeneous searchers whereas we extend the linearization approach to multiple classes of searchers, a camouflaging target, and explicitly include operational constraints about endurance and deconfliction.

At any time $t\in\cT$, the target moves according to a transition matrix $\Gamma_t$ whose element $\gamma_{s,c,s',c',t}$ represents the probability that a target occupying $(s,c)$ in period $t$ will be in $(s',c')$ during time period $t+1$. Contrary to {\bf CSP-U} and {\bf CSP-L}, the a priori enumeration of all possible target paths is not necessary in the following linearization. We adopt the additional notation in Table \ref{T4}.

\begin{table}[h]
	\centering
	\caption{Additional notation for model {\bf MSP}}
	\label{T4}
	\setlength\arrayrulewidth{1.25pt}
		\begin{tabular}{ll}
	  \hline
        \textbf{Indices}  \\
        \hline
        $j$ & Total search effort $j \in \{0, 1, 2, \dots, \}$ \\
        \hline
        \textbf{Sets}  \\
        \hline
        $\cJ_{s,t}^R$ & $= \{1, \dots, m_{s,t}\}$ \\
        \hline
        \textbf{Parameters}  \\
        \hline
        ${\alpha}_{c}$ & $\alpha_{c}= \alpha$ if $c = 0$ and $\alpha_{c}=0$ otherwise\\
        $\gamma_{s,c,s',c',t}$ & \multicolumn{1}{p{11.95cm}}{Probability that a target in state $(s,c)$ in period $t$ will be in state $(s',t')$ in period $t+1$} \\
        $p_{s,c}$  & \multicolumn{1}{p{11cm}}{Probability that the target is in state $(s,c)$ in period 1} \\
        $q_{s,c,t}$  & \multicolumn{1}{p{11cm}}{Probability that the target is in state $(s,c)$ in period $t$, i.e., $q_{s,c,t} = \sum_{s',c'} q_{s',c',t-1}\gamma_{s',c',s,c,t-1}$, $t=2, 3, \ldots, T$; $q_{s,c,1} = p_{s,c}$}\\
        $m_{s,t}$ & \multicolumn{1}{p{11cm}}{Maximum search effort possible in state $s$ at $t$: $m_{s,t} = \sum_{l \in \cL} m_{l,s,t}$}\\
        \hline
        \textbf{Decision Variables}  \\
        \hline
        $V_{s,t,j}$      & \multicolumn{1}{p{11.5cm}}{Binary variable = 1 if state $s$ receives $j$ search effort in period $t$, and = 0 otherwise}\\
        $P_{s,c,t}$        & Probability that target is in $(s,c)$ in $t$ and was not detected prior to $t$\\
        $Q_{s,c,t,j}$      & \multicolumn{1}{p{11cm}}{Auxiliary variable =  $P_{s,c,t}(1-e^{-j\alpha_c})$ if $V_{s,t,j}=1$ and = 0 otherwise}\\
        $W_{s,c,t}$        & \multicolumn{1}{p{11cm}}{Auxiliary variable = $P_{s,c,t}e^{-j{\alpha}_{c}}$ if $V_{s,t,j}=1$ and = $P_{s,c,t}$ otherwise} \\
        \hline
		\end{tabular}
\end{table}

We derive the linearization by introducing an ``information state" $P_{s,c,t}$ which represents the probability that the target occupies $(s,c)$ in period $t$ and that it has not been detected prior to $t$. We recall from {\bf SP} that $\sum_{l\in \cL} Z_{l,s,t}$ is the total search effort in state $s$ at period $t$. It is a nonnegative integer and can be represented equivalently by the binary variables $V_{s,t,j}$, each of which equals to 1 if there is $j$ search effort in state $s$ in period $t$, and equals to 0 otherwise.
This allows us to calculate the probability of detection over the entire time horizon as
\begin{equation}
\label{eqn:pd}
\sum_{t\in\cT} \sum_{(s,c) \in \mathcal{S} \times \{0,1\}}  P_{s,c,t}
\Bigg(1-\exp\Big(-\alpha_{c} \sum_{j \in \cJ^R_{s,t}}  j \ V_{s,t,j}\Big)\Bigg),
\end{equation}
where $\alpha_c = \alpha$ if $c=0$ and $\alpha_c = 0$ otherwise and $\cJ^R_{s,t} = \{1, \dots, m_{s,t}\}$, with $m_{s,t} = \sum_{l\in \cL} m_{l,s,t}$. The information state $P_{s,c,t}$ depends on the search plan as follows. The probability that the target occupies $(s,c)$ initially is $P_{s,c,1} = p_{s,c}$, which is an input parameter; see Table \ref{T4}. Moreover, it follows from the definition of $P_{s,c,t}$ and the Markov assumption that
\begin{equation}
\label{eqn:Pupdate}
P_{s,c,t+1} = \sum_{(s',c')\in \mathcal{S} \times \{0,1\}} \gamma_{s',c',s,c,t}P_{s',c',t}
\exp \Big(- \sum_{j\in \cJ^R_{s',t}} {\alpha}_{c'} \ j \ V_{s',t,j}\Big)
\end{equation}
for $s,c$ and $t=1, 2, \ldots, T-1$.

We shall linearize the nonlinear expressions \eqref{eqn:pd} and \eqref{eqn:Pupdate}.
First, we linearize the probability of non-detection (i.e., the complement of \eqref{eqn:pd}) via the introduction of the auxiliary variable $Q_{s,c,t,j}$ which takes value
$P_{s,c,t}(1-e^{-j {\alpha}_{c}})$ if $V_{s,t,j}=1$ and takes value 0 otherwise.
This linearization is accomplished using constraints
\eqref{splm:Zlin12} and \eqref{splm:Zlin22} below.
The inequality \eqref{splm:Zlin12} is a ``big-M'' constraint where any constant at least as large as $P_{s,c,t}$ is needed to multiply $(1-e^{-j {\alpha}_{c}})$.
Since $P_{s,c,t}$ is the probability that the target is in $(s,c)$ in period $t$ {\it and} that the target is not detected prior to $t$ and
$q_{s,c,t}$ is the probability that the target is in $(s,c)$ in period $t$ as defined in Table \ref{T4}, we must have $q_{s,c,t}\geq P_{s,c,t}$ for all $(s,c) \in \cS \times \{0,1\}, t \in \cT$.
Consequently, each ``big-M" parameter in \eqref{splm:Zlin12} is set to
$q_{s,c,t}$. Using the same rationale, we let $q_{s,c,t}$ furnish the bound on $P_{s,c,t}$ in \eqref{splm:PUB2} below.
Second, the evolution of the information state is also nonlinear as it can be seen from \eqref{eqn:Pupdate}.
We linearize that expression by means of the auxiliary variable $W_{s,c,t}$ and constraints
\eqref{splm:PW2}-\eqref{splm:Wlin2} below.
Note that $W_{s,c,t}$ is equal to
$P_{s,c,t}e^{-j {\alpha}_{c}}$ if $V_{s,t,j}=1$ and is equal to  $P_{s,c,t}$ otherwise. Compiling these derivations, we obtain the following equivalent MILP reformulation of {\bf SP} under the Markovian target path model.


\begin{subequations} 
\begin{align}
\textbf{MSP:}~~~~~&\nonumber\\
\nnmin_{X,P,Q,V,W} & \; 1- \sum_{(s,c) \in \cS \times \{0,1\} } \sum_{t \in \cT} \sum_{j \in \cJ^R_{s,t}} Q_{s,c,t,j} \label{splm:SP1LMobj2}\\
\text{subject to } & \; Q_{s,c,t,j} \leq q_{s,c,t}(1-e^{-j{\alpha}_{c}}) V_{s,t,j} & (s,c) \in \cS \times \{0,1\},  t \in \cT, j\in \cJ^R_{s,t}\label{splm:Zlin12}\\
& Q_{s,c,t,j} \leq (1-e^{-j{\alpha}_{c}})P_{s,c,t}  & (s,c) \in \cS \times \{0,1\},  t \in \mathcal{T}, j\in \cJ^R_{s,t} \label{splm:Zlin22}\\
& P_{s,c,t+1} = \sum_{(s',c') \in \cS\times \{0,1\}} \gamma_{s',c',s,c,t}W_{s',c',t}  & (s,c) \in \cS \times \{0,1\},  t \in \cT \setminus \{T\} \label{splm:PW2}\\
& W_{s,c,t} \leq P_{s,c,t}  & (s,c) \in \cS \times \{0,1\},  t \in \cT \label{splm:WUB2}\\
& W_{s,c,t} \leq e^{-j{\alpha}_{c}}P_{s,c,t} + q_{s,c,t}(1-e^{-j{\alpha}_{c}})(1-V_{s,t,j})  &  (s,c) \in \cS \times \{0,1\},  t \in \mathcal{T}, j\in \cJ^R_{s,t} \label{splm:Wlin2}\\
& P_{s,c,1} = p_{s,c} &  (s,c) \in \cS \times \{0,1\} \label{splm:PUB12}\\
& P_{s,c,t}\leq q_{s,c,t} & (s,c) \in \cS \times \{0,1\}, \  t \in \cT \label{splm:PUB2}\\
& \sum_{l \in \cL} \sum_{\substack{s'\in\cR_l(s) \\ t-d_{l,s',s}\geq 0}} \beta_{l,s',s,t}X_{l,s',s,t-d_{s',s}} = \sum_{j\in \cJ^R_{s,t}} j \ V_{s,t,j} & {s \in \cS, t \in\cT} \label{splm:SPI1LMall2}\\
& \sum_{j\in \cJ^R_{s,t}} V_{s,t,j} = 1 & {s \in \cS, t \in \cT} \label{splm:SPI1LMall22}\\
& (\ref{eqn:SPXflow})\mbox{-}(\ref{eqn:takeoff2}); \eqref{NEW2}\mbox{-}\eqref{REL-INTEb} \; \; \nonumber\\
& P_{s_c,t},W_{s_c,t} \geq 0 & (s,c) \in \cS \times \{0,1\}, t \in \cT\\
& Q_{s_c,t,j} \geq 0 & (s,c) \in \cS \times \{0,1\}, t \in \cT, j\in \cJ^R_{s,t} \\
& V_{s,t,j}\in\{0, 1\} & s \in \cS, t \in \cT, j\in \cJ^R_{s,t}
\end{align}
\end{subequations}
The objective function \eqref{splm:SP1LMobj2} gives the probability of non-detection; its correctness follows from \eqref{eqn:pd}. The binary variable $V_{s,t,j}$ is linked to $X_{l,s,s',t}$ in \eqref{splm:SPI1LMall2}. The remaining constraints follow from the discussion above.\\

\state Computational Tests. We consider two instances of {\bf MSP} of the kind described in Subsection  \ref{RF1}, but now with the Markovian target path model obtained from the transition probabilities described there. This produces the last row of Table \ref{tableapro2} for the two instances that only differ in the number of searchers ($J_1$) and the planning horizon ($T$). Neither instance of {\bf MSP} can be solved directly using Gurobi within 900 seconds. While an optimal solution is eventually achieved in the instance with $J_1=3$, $T=12$, the gap is sizable in the other instance after 900 seconds; the lower bound is 0.4043 and the upper bound 0.4659 at that time. We conclude that {\bf MSP} is computationally challenging and this motivates the derivation of cutting plane algorithms in the next subsection.

Table \ref{tableapro2} also illustrates how the Markovian target path model can be viewed as the limit of the conditional target path models when the latter are obtained by sampling according to the Markov transition matrices. With a planning horizon of $T=12$ and the present Markovian target path model with typically 5 possible moves per time period, we obtain that the model produces about $|\Omega| = 5^{12} \approx 2 \cdot 10^8$ target paths. Thus, the sample sizes ranging from 100 to 5000 in Table \ref{tableapro2} are relatively small. Nevertheless, the sample average approximations have minimum objective function values close to those for the Markovian target path model when the sample size is at least 1000. (This motivates in part our focus on conditional target path models with 1000 paths in Section \ref{SGM1}.) There is a significant computational advantage of considering sample averages; Section \ref{SGM1} provides extensive evidence that conditional target path models are tractable. Table \ref{tableapro2} provides a direct comparison using {\bf CSP-L-Pre} as the approach for solving the sample average approximations. Further speed-up might be possible with {\bf CSP-U-Pre} or the {\tt OA} method.

\begin{table}[H]
	\centering
	\caption{{\small For $S=82$ states and varying numbers of sampled target paths: Min-value and solution time (sec.) to relative optimality gap of 0.0001 or, if not reached in 900 seconds, relative optimality gap in brackets after 900 seconds. The case marked with asterisk solves to optimality in 1604 seconds.}}
	{\small\begin{tabular}{|c|r|r|r|r|}
		\hline
		& \multicolumn{2}{c|}{$J_1=3$, $T=12$} &  \multicolumn{2}{c|}{$J_1=5$, $T=10$}\\
		\hline
		Sample size  & Min-value & Solution time & Min-value & Solution time\\
		\hline
		100   & 0.2931& 0.6       & 0.3039   & 0.4     \\
		500   & 0.4048& 0.5         & 0.4032   & 2      \\
		1000  & 0.5007& 2          & 0.4180   & 22     \\
		2000  & 0.5031& 6         & 0.4336  & 68        \\
		5000  & 0.4973 & 246       & 0.4266   & 16     \\
		\hline
		Markovian & 0.5036 & *[0.0332]  & 0.4043-0.4659 & [0.0916] \\
		\hline
	\end{tabular}}
	\label{tableapro2}
\end{table}

\subsection{Cutting Plane Algorithms} \label{AL2}

In this subsection, we extend the cutting plane methods of \cite{Royset2010} to the present setting with a camouflaging target and heterogenous searchers. A direct extension yields {\tt SCA} in Subsection \ref{A1}. Further refinements leveraging bundles and outer approximations follow in Subsections \ref{A2} and \ref{A3}.

\subsubsection{Secant Cutting Plane Algorithm (SCA)}\label{A1}

Adaptively constructed piecewise-linear approximations of the objective function in {\bf SP} lead to a cutting plane method {\tt SCA}, which in each iteration $i$ solves the MILP:
\begin{subequations}
\begin{align}
{\bf P^i_{SCA}}: \; \;
& \nnmin \; \xi \notag\\
\text{subject to } & \ \xi \geq f(Z^k)
+ \sum\limits_{l\in \cL} \sum\limits_{s \in \cS} \sum\limits_{t \in \cT}
(f(Z^k + \Delta_{l,s,t}) - f(Z^k))(Z_{l,s,t} - Z^k_{l,s,t}), ~~  k=1,\ldots,i \label{sca:Pi2}\\
&\; \eqref{eqn:SPXflow}\mbox{-}\eqref{NEW3}\notag
\end{align}
\end{subequations}
where $f(Z)$ denotes the objective function of {\bf SP} and $Z^k$ is the allocation of search effort from a previous iteration. The notation $\Delta_{l,s,t} \in \{0,1\}^{\cL\times \cS\times\cT}$ refers to a Boolean parameter vector in which all elements are 0 except the $(l,s,t)$-component equals to 1 and is used to measure the impact of varying one single variable $Z_{l,s,t}$ on the value of the objective function.
A new secant cut \eqref{sca:Pi2} is added at each iteration $i$ and problem ${\bf P^i_{SCA}}$ minimizes the resulting piecewise-linear approximation of $f(Z)$.

Guided by \cite{Royset2010}, the calculation of a secant cut proceeds in two steps: (i) compute the probability $r_{s,c,t}(Z)$ that the target is in $(s,c)$ at time $t$ and is not detected before $t$, and (ii) compute the probability $\bar{r}_{s,c,t}(Z)$ that the target is not detected in the periods after $t$ given that the target is in $(s,c)$ at time $t$.
We define $r_{s,c,1}(Z)=p_{s,c}$ and $\bar{r}_{s,c,T}(Z)=1$ so that all other $r_{s,c,t}(Z)$ and $\bar{r}_{s,c,t}(Z)$ can be calculated recursively as follows:
\begin{subequations}
\begin{align}
r_{s,c,t}(Z) & = \sum_{s',c'} r_{s',c',t-1}(Z)~ \gamma_{s',c',s,c,t-1}~ e^{-\sum_{l \in \cL} \alpha_{c'}
\label{CUT1}
Z_{l,s',t-1}} \\
\bar{r}_{s,c,t}(Z) & = \sum_{s',c'} \bar{r}_{s',c',t+1}(Z)~ \gamma_{s,c,s',c',t}~ e^{-\sum_{l \in \cL} \alpha_{c'} \label{CUT2}
Z_{l,s',t+1}}.
\end{align}
\end{subequations}
This allows us in turn to calculate, for any $t\in \cT$, the objective function
\begin{equation}
\label{FT}
f(Z)= \sum_{s,c} r_{s,c,t}(Z)~ e^{-\sum_{l \in \cL} \alpha_c 
Z_{l,s,t}}~ \bar{r}_{s,c,t}(Z),
\end{equation}
which is the product of the probability of not being detected before $t$, the non-detection probability at $t$, and the probability of not being detected after $t$.
A secant cut can then be computed via
\begin{equation*}
	\label{FTD}
	f(Z+\Delta_{l,s,t})-f(Z) = r_{s,0,t}(Z)
		\bigg( e^{-\sum_{l \in \cL} \alpha (Z_{l,s,t}+1)}
		-e^{-\sum_{l \in \cL} \alpha Z_{l,s,t}} \bigg) \ \bar{r}_{s,0,t}(Z) \ .
\end{equation*}
This derivation deviates from that of \cite{Royset2010} by accounting for a camouflaging target and heterogeneous searchers.

We can now present the formal structure of {\tt SCA}. Let $\delta, \delta_i \geq 0, i=0,1,2,\dots,N$ denote optimality tolerances while $\underline{\xi}$ and $\overline{\xi}$ are lower and the upper bounds on the optimal value of {\bf SP}.

\vspace{0.1in}
\noindent
{\bf \underline{Initialization:}}
\newline
\textbf{Step 0}: Set: $\underline{\xi}=0$; $\overline{\xi}=1$; $i = 1$; $Z^1 = 0$ (zero vector). \\
{\bf \underline{Iterative process -- Iteration $i$:}} \\
\textbf{Step 1}:
Calculate $f(Z^i)$. If $f(Z^i)< \overline{\xi}$, then set $\overline{\xi}= f(Z^i)$. \\
\textbf{Step 2}:
If $\overline{\xi} - \underline{\xi} \leq \delta \underline{\xi}$, then stop: tolerance satisfied. \\
\textbf{Step 3}:
Solve problem ${\bf P^i_{SCA}}$ to tolerance $\delta_i$, achieve solution $Z^{i+1}$, and lower bound $\underline{\xi}^{i+1}$. \\
\textbf{Step 4}:
If $\underline{\xi}^{i+1} > \underline{\xi}$, then set $\underline{\xi} = \underline{\xi}^{i+1}$. \\
\textbf{Step 5}:
If $\overline{\xi} - \underline{\xi} \leq \delta \underline{\xi}$, then stop: tolerance satisfied.
Else, replace $i$ with $i+1$ and go to Step 1.

\vspace{0.1in}

In the numerical tests of Subsection \ref{NT2}, we set $\delta_1 = 0$ and $\delta_i = \min\{0.03, g_i/3\}$ for $i \geq 2$, where $g_i = (\overline{\xi} - \underline{\xi})/\underline{\xi}$ is computed after Step 1 of
iteration $i$. However, we use $\delta_i = \min\{0.03, g_i/3, \delta_{i-1}/2\}$ if
$Z_i$ is a repetition of a previously obtained solution.

\subsubsection{Bundle-based Cutting Plane Algorithm (B-SCA)}\label{A2}

We refine {\tt SCA} by incorporating bundles as well as preprocessing techniques. As a preliminary step, we partition the set $\cT$ into two mutually exclusive subsets
\begin{equation}
\label{SET2}
\mathcal T_{nd} = \left\{t \in \cT:
\sum\limits_{s \in \mathcal{S}} q_{s,0,t} \ \beta_{s,t} =0 \right\}
\end{equation}
that includes all periods (and only those) at which no detection can occur, and its complement $\cT_d = \cT \setminus \cT_{nd}$ which includes the periods at which detection is possible.
The notation $\beta_{s,t}$ in \eqref{SET2} specifies a Boolean parameter with value 0 if no searchers can reach state $s$ by time $t$ and value 1 otherwise. For each $t \in \cT_d$, we build the set
\begin{equation*}
\label{ST}
\mathcal{V}^{t}_{nd} = \left\{s \in \cS:
q_{s,0,t} \ \beta_{s,t} = 0 \right\} \ , \ t \in \cT_{d}
\end{equation*}
that contains all states $s$ for which no detection can occur at $t$. We use the notation $\mathcal{V}^{t}_{d}$ to refer to the complement of $\mathcal{V}^{t}_{nd}$:   $\mathcal{V}^{t}_{d} = \cS\ \setminus \mathcal{V}^{t}_{nd}$.

The above sets are used via a bundling approach to reduce the size of the decision and constraint spaces.
First, we eliminate the integer decision variables $Z_{l,s,t}$ at any period $t \in \cT_{nd}$ when no detection can occur across all states. Since no detection can occur at these periods, we do not need to keep track of how many searchers are in $s$ at $t$.
Second, at the {\it remaining} periods $t \in \cT_d$, we further remove the integer decision variable $Z_{l,s,t}$ for any $(s,t) \in \mathcal{V}^{t}_{nd}, t \in \cT_{d}$ corresponding to any state $s$ at which detection is impossible.
More precisely,  for any $t \in \cT_d$, we combine all tuples $(s,t), s \in \mathcal{V}^{t}_{nd}, l \in \cL$ in a so-called {\it bundle} $\cB^t$ and do not include any variable $Z_{l,s,t}$ for any tuple $(s,t)$ included in one of the bundles  $\cB^t, t \in \cT_d$. This produces the algorithm {\tt B-SCA}, which in each iteration $i$ solves the MILP:
\begin{subequations}
\begin{align}
& {\bf P^i_{B-SCA}}: \; \; \nnmin \; \xi \notag\\
\text{subject to } & \ \xi \geq f(Z^k) + \sum_{\substack{t \in \cT_{d}, l\in \cL, \\ s \in \mathcal{V}^{t}_{d}}}
(f(Z^k + \Delta_{l,s,t}) - f(Z^k))(Z_{l,s,t} - Z^k_{l,s,t}), ~k=1,\ldots,i \label{eqn:Pi2} \\
&\sum_{\substack{s'\in\cR_l(s) \\ t-d_{l,s',s}\geq 0}} \beta_{l,s',s,t} X_{l,s',s,t-d_{l,s',s}} = Z_{l,s,t}, \qquad t \in \cT_{d}, s \in \mathcal{V}^{t}_{d}, l\in \cL \label{eqn:Pi5}\\
&\eqref{eqn:SPXflow}\mbox{-}\eqref{eqn:takeoff2}; \eqref{NEW2}\mbox{-}\eqref{REL-INTEb}\notag \\
& Z_{l,s,t} \in \{0, 1, 2, \ldots, m_{l,s,t}\}  , \qquad  t \in \cT_d, s \in \mathcal{V}^{t}_{d}, l\in \cL. \label{eqn:Pi9}
\end{align}
\end{subequations}

Due to the smaller number of variables $Z_{l,s,t}, s \in \cV^t_{d}, t \in \cT_d$ used by {\tt B-SCA}, we can simplify \eqref{CUT1}-\eqref{CUT2} as follows:
\begin{align}
\hspace{-0.09in}	
r_{s,c,t}(Z) = \left\lbrace
\begin{matrix*}[l]
0 \; & \; \text{if} \  q_{s,0,t} = 0  \\
\sum\limits_{s',c'} r_{s',c',t-1}(Z)~ \gamma_{s,c,s',c',t-1}  \; & \; \text{if} \ \beta_{s,t} = 0 \ \text{and} \ q_{s,0,t} \neq 0  \\
\sum\limits_{s',c'} r_{s',c',t-1}(Z)~ \gamma_{s,c,s',c',t-1}
e^{-\sum_{l \in \cL} \alpha_{c'} Z_{l,s', t-1}} \; & \; \text{otherwise}
\end{matrix*}
\right. \notag
\end{align}
\begin{align}
\hspace{-0.15in}		
\bar{r}_{s,c,t}(Z) = \left\lbrace
\begin{matrix*}[l]
0 \; & \; \text{if} \ q_{s,0,t} = 0  \\
\sum\limits_{s',c' } \bar{r}_{s',c',t+1}(Z)~ \gamma_{s_c,s'_c,t}  \; & \; \text{if} \
 \beta_{s,t} = 0 \ \text{and} \ q_{s,0,t} \neq 0  \\
\sum\limits_{s',c'} \bar{r}_{s',c',t+1}(Z)~ \gamma_{s,c,s',c',t}~ e^{-\sum_{l \in \cL} \alpha_{c'}
Z_{l,s',t+1}} \; & \; \text{otherwise.}
   \end{matrix*}
     \right. \notag
\end{align}

\subsubsection{Bundle-based Cutting Plane Algorithm with Outer Approximation (OA-B-SCA)}\label{A3}

We next adjust {\tt B-SCA} by replacing the feasible sets of each subproblem by an outer approximation. While the outer approximation remains mixed-integer, it can be described by fewer integer variables and constraints.
The expectation is that the size reduction of the decision and constraint spaces will allow for a quicker solution of the subproblems. The trade-off is that the feasible sets of the subproblems are relaxations and will therefore provide looser lower bounds on the optimal value of the actual problem.

The resulting algorithm {\tt OA-B-SCA} rests on the following rationale.
We  observe that the probability $q_{s,0,t}$ of a target being in $(s,0)$ at $t$ can significantly vary across pairs $(s,t)$. Even if positive, some $q_{s,0,t}$ can be extremely low making it ineffective to place a searcher in $s$ at $t$. Building on this, each subproblem in the proposed outer-approximation algorithm {\tt OA-B-SCA} leverages integer decision variables $Z_{l,s,t}$ only for tuples $(s,t)$ with the largest $q_{s,0,t}$ across all states $s$ at $t$, i.e., the states where a target is most likely to be at $t$.
As for {\tt B-SCA}, we first drop the integer variables $Z_{l,s,t}$ for any tuple $(l,s,t)$ with $t \in \cT_{nd}$.
We then remove the integer variables $Z_{l,s,t}$ corresponding to the tuples $(l,s,t)$ for any $(s,t)$ pairs at which detection is impossible and those at which probability of the target being in state $s$ at time $t$ is not one of the highest.

We denote by $\mathcal W_{t,\upsilon}$ the set of tuples $(s,t)$ associated with the $\upsilon$ most likely states for the target to be in and not be camouflaging at time $t$.
Let $\cW^c_{t,\upsilon}$ be its complement.
For each $(s,t) \in \cW^c_{t,\upsilon}, l \in \cL$, we relax the integrality condition on the variables $Z_{l,s,t}$. This produces the algorithm {\tt OA-B-SCA}, which in each iteration $i$ solves the subproblem:
\begin{subequations} \notag
\begin{align}
{\bf P^i_{OA-B-SCA}}: \
\nnmin \; & \xi \notag \\
    \text{subject to } \; &  \eqref{eqn:SPXflow}\mbox{-}\eqref{eqn:takeoff2}; \eqref{NEW2}\mbox{-}\eqref{REL-INTEb} ; \eqref{eqn:Pi2}\mbox{-}\eqref{eqn:Pi5} \notag \\
    & Z_{l,s,t} \in \{0, 1, 2,  \ldots, m_{l,s,t}\}  , \qquad  t \in \cT_d, (s,t) \in \mathcal{W}_{t, \upsilon}, l\in \cL \notag\\
    & Z_{l,s,t} \in [0, m_{l,s,t}]  , \qquad \qquad  \quad t \in \cT_d, (s,t) \in \mathcal{W}^{c}_{t, \upsilon}, l\in \cL. \notag
\end{align}
\end{subequations}
The feasible set of each subproblem ${\bf P^i_{OA-B-SCA}}$ is a {\it relaxation} of the actual feasible set of {\bf SP}. As with {\tt SCA} and {\tt B-SCA}, the feasible set of {\tt OA-B-SCA} is defined by mixed-integer linear constraints, but it contains  (many) {\it fewer} integer variables than the feasible sets of {\tt SCA} and {\tt B-SCA}.

The structure of {\tt OA-B-SCA} is similar to that of {\tt B-SCA}. However, the stopping criterion differs. Due to the relaxation of the integrality restrictions of a subset of the variables $Z_{l,s,t}$, the solution obtained by solving the subproblems ${\bf P^i_{OA-B-SCA}}$ is not necessarily feasible for {\bf MSP} and a post-optimization step must be carried out to restore feasibility and allow for the computation of a valid upper bound.

If the solution of ${\bf P^i_{OA-B-SCA}}$ is fractional, we do not have a valid upper bound. To obtain one, we must first restore integrality, which can be done in a heuristic manner, by using a basic rounding procedure, or by solving a reduced-size {\it integrality restoration} problem.
The integrality restoration problem is a much simplified variant of ${\bf P^i_{OA-B-SCA}}$ and contains many less integer variables so that it can be solved to optimality extremely quickly (typically in less than one second). Actually, we do not need to solve it to optimality since any feasible solution provides a valid upper bound for the true problem.

Let $\bar{Z}^i$ be the solution produced by ${\bf P^i_{OA-B-SCA}}$ at iteration $i$.
We fix all variables $Z_{l,s,t}$ which have an integer value in $\bar{Z}^i_{l,s,t}$ and they become fixed parameters. 
Denoting by $\mathbb Z_+$ the set of nonnegative integers, we define
\[
\mathcal{A}^{I}_i = \big\{(l,s,t)\in \cL \times \cV^t_d \times \cT_d: \bar{Z}^i_{l,s,t} \in \mathbb Z_+\big\}, ~~ \mathcal{A}^{F}_i  =  \big\{(l,s,t)\in \cL \times \cV^t_d \times \cT_d: \bar{Z}^i_{l,s,t} \notin \mathbb Z_+\big\}.
\]
The sets $\cA^I_i$ and $\cA^F_i$ include the tuples $(l,s,t)$ whose corresponding variables $Z_{l,s,t}$ respectively take integer and fractional values $\bar{Z}^i_{l,s,t}$ in the obtained solution of ${\bf P^i_{OA-B-SCA}}$. The sets $\cA^I_i$ and $\cA^F_i$ are updated at each iteration $i$.
The reduced-size MILP integrality restoration subproblem ${\bf IR^i}$ at $i$ then reads:
\begin{subequations} \notag
\begin{align}
{\bf IR^i}: \
\nnmin \; & \xi \notag \\
    \text{subject to } \; &  \eqref{eqn:SPXflow}\mbox{-}\eqref{eqn:takeoff2}; \eqref{NEW2}\mbox{-}\eqref{REL-INTEb} ; \eqref{eqn:Pi2}\mbox{-}\eqref{eqn:Pi5} & \notag \\
    & Z_{l,s,t} = \bar{Z}^i_{l,s,t}   & (l,s,t) \in \cA^I_i \notag\\
    & Z_{l,s,t} \in \{0, 1, 2, \ldots, m_{l,s,t}\}  & (l,s,t) \in \cA^F_i. \notag
\end{align}
\end{subequations}

\noindent
The algorithm {\tt OA-B-SCA} is structured as follows:

\noindent
{\bf \underline{Initialization:}}
\newline
\textbf{Step 0}: Set: $\underline{\xi}=0$; $\overline{\xi}=1$; $i = 1$; $Z^1 = 0$ (zero vector). \\
{\bf \underline{Iterative process -- Iteration $i$:}} \\
\textbf{Step 1}:
Calculate $f(Z^i)$. If $f(Z^i)< \overline{\xi}$, then set $\overline{\xi}= f(Z^i)$. \\
\textbf{Step 2}:
If $\overline{\xi} - \underline{\xi} \leq \delta \underline{\xi}$, then stop: tolerance satisfied. \\
\textbf{Step 3}:
Solve problem ${\bf P^i_{OA-B-SCA}}$  to tolerance $\delta_i$, achieve solution $\bar Z^{i}$, and lower bound $\underline{\xi}^{i+1}$ \\
\textbf{Step 4}:
If $\underline{\xi}^{i+1} > \underline{\xi}$, then set $\underline{\xi} = \underline{\xi}^{i+1}$. \\
\textbf{Step 5}:
If $\overline{\xi} - \underline{\xi} \leq \delta \underline{\xi}$, then stop:  tolerance satisfied.\\
\textbf{Step 6}: If $\bar{Z}^i$ is integer, set $Z^{i+1} = \bar Z^i$. Else, solve ${\bf IR^i}$ to restore integrality and obtain $Z^{i+1}$.\\
\textbf{Step 7}:
Replace $i$ with $i+1$ and go to Step 1.

\subsection{Numerical Tests} \label{NT2}

We compare the three cutting plane methods {\tt SCA}, {\tt B-SCA}, and {\tt OA-B-SCA} with a direct solution of {\bf MSP} across two groups of instances.\\

\state Homogenous Searchers. We first consider problem instances of the kind described in Subsection \ref{RF1}, except that we consider here a Markovian target path model. These instances do not allow for the camouflage option and there is no endurance limit. Table \ref{table7} reports the computational time for three searchers, 82 states, and varying planning horizon $T$. For instances with few time periods ($T \leq 11$), the direct solution of \textbf{MSP} is faster than the cutting plane methods {\tt SCA}, {\tt B-SCA}, and {\tt OA-B-SCA}. As $T$ increases beyond 11, the optimality gap with the three cutting plane methods is smaller. In particular, for all instances with 12 or more periods, the outer-approximation algorithm {\tt OA-B-SCA} performs best and reduces the optimality gap the most. For $T$ = 13 (resp., 14 and 15),
{\tt OA-B-SCA} produces a gap of 0.0161 (resp., 0.0239 and 0.0183) less than {\tt SCA}.
 These results highlight the efficiency of {\tt OA-B-SCA} in solving the most challenging instances of this type.

\begin{table}[H]
\centering
\caption{{\small For Markovian target model, $S=82$ states, $J_1=3$ searchers, and varying numbers of time periods: Solution time (sec.) to relative optimality gap of 0.0001 or, if not reached in 900 seconds, relative optimality gap in brackets after 900 seconds.}}
    {\small\begin{tabular}{|c|r|r|r|r|}
     \hline
     $T$ & {\bf MSP} & {\tt SCA} & {\tt B-SCA} & {\tt OA-B-SCA} \\
     \hline
     7 & 0.1       & 5        & 5        & 5 \\
     8 & 0.3       & 46       & 37       & 37 \\
     9 & 0.8       & 87       & 66       & 64 \\
     10& 9         & [0.0186] & [0.0175] & [0.0198] \\
     11& 278       & [0.0590] & [0.0574] & [0.0581] \\
     12& [0.1693]  & [0.0983] & [0.1006] & [0.0916] \\
     13& [0.3151]  & [0.1410] & [0.1316] & [0.1249]    \\
     14& [0.4257]  & [0.1742] & [0.1742] & [0.1503] \\
     15& [0.5357]  & [0.1915] & [0.1969] & [0.1732] \\
     \hline
     Average Optimality Gap & 0.1606 & 0.0758 & 0.0753 & 0.0686 \\
     \hline
\end{tabular}}
\label{table7}
\end{table}

Table \ref{table8} considers instances with $J_1=15$ searchers. As observed in Table \ref{table7}, solving \textbf{MSP} directly is the most computationally efficient approach for small instances ($T = 7$ and possible 8) but the three cutting plane algorithms dominate \textbf{MSP} when the planning horizon increases and the instances become more challenging. Among the three, {\tt B-SCA} is the most efficient on most instances, but is closely followed by {\tt OA-B-SCA}.
On average, for the challenging instances ($T \geq 9$), the optimality gap with {\tt B-SCA} is on average 0.0022 lower than for {\tt SCA}, which highlights the computational benefits of the bundle-based cutting plane {\tt B-SCA}.

\begin{table}[H]
\centering
\caption{{\small For Markovian target model, $S=82$ states, $J_1=15$ searchers, and varying numbers of time periods: Solution time (sec.) to relative optimality gap of 0.0001 or, if not reached in 900 seconds, relative optimality gap in brackets after 900 seconds.}}
    {\small\begin{tabular}{|c|r|r|r|r|}
     \hline
     $T$ & {\bf MSP} & {\tt SCA} & {\tt B-SCA} & {\tt OA-B-SCA} \\
     \hline
     7  & 2         & 8        & 7        & 8 \\
     8  & 95        & 204      & 83       & 75 \\
     9  & [0.0624]  & [0.0015] & [0.0005] & [0.0007] \\
     10 & [0.2039]  & [0.0035] & [0.0032] & [0.0032] \\
     11 & [0.3502]  & [0.0054] & [0.0048] & [0.0047] \\
     12 & [0.5144]  & [0.0092] & [0.0078] & [0.0065] \\
     13 & [0.7010]  & [0.0146] & [0.0135] & [0.0203]    \\
     14 & [0.8783]  & [0.0259] & [0.0220] & [0.0258] \\
     15 & [1.1006]  & [0.0443] & [0.0332] & [0.0377] \\
     \hline
     Average Optimality Gap & 0.4211 & 0.0116 & 0.0094 & 0.0109 \\
     \hline
\end{tabular}}
\label{table8}
\end{table}

Table \ref{table9} examines the effect of the number of searchers on the solution time. For $J_1 \leq 4$, the direct solution of {\bf MSP} dominates the cutting plane approaches. However, {\tt SCA}, {\tt B-SCA}, and {\tt OA-B-SCA} have a clear advantage when the number of searchers exceeds 4. The algorithm {\tt B-SCA} is the best of the three on all instances, but the differences are modest.

\begin{table}[H]
\centering
\caption{{\small For Markovian target model, $S=82$ states, $T=10$ time periods, and varying numbers of searchers: Solution time (sec.) to relative optimality gap of 0.0001 or, if not reached in 900 seconds, relative optimality gap in brackets after 900 seconds.}}
    {\small\begin{tabular}{|c|r|r|r|r|}
     \hline
     $J_1$ & {\bf MSP} & {\tt SCA} & {\tt B-SCA} & {\tt OA-B-SCA} \\
     \hline
     1 & 0.3       & 34       & 34       & [0.0363]  \\
     2 & 1         & [0.0017] & 581      & [0.0159]   \\
     3 & 9         & [0.0186] & [0.0175] & [0.0213]   \\
     4 & 70        & [0.0236] & [0.0213] & [0.0245]   \\
     5 & [0.0340]  & [0.0163] & [0.0161] & [0.0184]  \\
     10& [0.1400]  & [0.0060] & [0.0052] & [0.0057]  \\
     15& [0.2000]  & [0.0035] & [0.0032] & [0.0032] \\ \hline
     Average Optimality Gap& 0.0534 & 0.0099 & 0.0090 & 0.0179 \\
     \hline
\end{tabular}}
\label{table9}
\end{table}

Table \ref{table10} examines the effect of the size of the square grid of cells and thus the number of states. For small grid sizes (i.e., less than 7-by-7 cells producing $S \leq 50$), the cutting-plan approaches dominate the direct solution of \textbf{MSP}. However, the direct solution of \textbf{MSP} is by far the fastest approach to prove optimality for larger grid sizes, such as $S \geq 82$, which turns out to be the simplest instances. The approach solves all those instances with an average solution time of 2.3 seconds whereas the three cutting plane methods struggle to solve the 82-state instance and are slower to prove optimality for the three instances with $S$ = 122, 170, and 226.
Among the cutting plane methods, {\tt OA-B-SCA} has the lowest average optimality gap when optimality cannot be proven and has the smallest average solution time for the other instances.

\begin{table}[h]
\centering
\caption{{\small For Markovian target model, $J=3$ searchers, $T=10$ time periods, and varying numbers of states: Solution time (sec.) to relative optimality gap of 0.0001 or, if not reached in 900 seconds, relative optimality gap in brackets after 900 seconds.}}
    {\small\begin{tabular}{|c|r|r|r|r|}
     \hline
     $S$ & {\bf MSP} & {\tt SCA} & {\tt B-SCA} & {\tt OA-B-SCA} \\
     \hline
     26 & [0.8314]  & [0.2355] & [0.2372] & [0.2290]  \\
     50 & [0.1521]  & [0.1070] & [0.1041] & [0.0995]   \\
     82 & 8         & [0.0186] & [0.0174] & [0.0168]   \\
     122& 0.6       & 86       & 80       & 75   \\
     170& 0.4       & 24       & 28       & 30  \\
     226& 0.3       & 15       & 18       & 16  \\
     \hline
     Average Optimality Gap & 0.1633  & 0.0602  & 0.0598 & 0.0575 \\
     \hline
\end{tabular}}
\label{table10}
\end{table}

To sum up, the results reported in Tables \ref{table7}-\ref{table10} demonstrate that while the linear reformulation \textbf{MSP} tends to be quicker for the smallest and least challenging instances, the two proposed bundle-based cutting plane algorithms {\tt B-SCA} and {\tt OA-B-SCA} are superior for the challenging ones. They also improve on {\tt SCA}, which in the present setting with homogenous searchers, no endurance constraints, and no camouflaging is essentially equivalent to an algorithm from \cite{Royset2010}.
It appears that, depending on the type of instances, it is preferable to derive stronger lower bounds (as allowed by {\tt B-SCA)} while, for others, a quicker solution time of the subproblems (as allowed by {\tt OA-B-SCA}) and thus the execution of more iterations within a given allowed time is more beneficial.\\

\state Heterogeneous Searchers and Camouflaging. We next consider instances of the kind associated with Table \ref{table:CamoEndurance}, which involves camouflaging, endurance constraints, and two classes of searchers. Table \ref{table-2class} presents the results for instances with $J=J_1+J_2=3$ and $J=J_1+J_2=15$ searchers. When $J=3$, we consider two searchers of class 1 and one searcher of class 2.  When $J=15$, we consider ten searchers of class 1 and five of class 2. The classes only differ in terms of endurance.

\begin{table}[H]
	\centering
	\caption{{\small For Markovian target model, $S=83$ states, and varying numbers of time periods and searchers across two classes with varying endurance: Solution time (sec.) to relative optimality gap of 0.0001 or, if not reached in 900 seconds, relative optimality gap in brackets after 900 seconds.}}
	{\small\begin{tabular}{|c|c|r|r|r|r|}
		\hline
		$J$ & $T$ & {\bf MSP} & {\tt SCA} & {\tt B-SCA} & {\tt OA-B-SCA} \\
		\hline
		3 & 10 & 4         & 58        & 73        & 64  \\
		3 & 12 & 78        & [0.0111]  & [0.0114]  & [0.0200] \\
		3 & 14 & [0.0953]  & [0.0630]  & [0.0696]  & [0.0219]  \\
		3 & 15 & [0.2510]  & [0.0814]  & [0.0809]  & [0.0557]  \\
		3 & 16 & [0.2183]  & [0.0726]  & [0.0655]  & [0.0209]  \\
		3 & 17 & [0.4995]  & [0.1222]  & [0.1470]  & [0.0360]  \\
		3 & 18 & [0.6250]  & [0.1292]  & [0.1483]  & [0.0330] \\
		3 & 20 & [2.6835]  & [0.1992]  & [0.2426]  & [0.1681] \\
		\hline
		& Average & 0.5466 &  0.0848 & 0.0957 & 0.0444 \\
		\hline
		15 & 10 & [0.0977]  &  22        &  23        & 13   \\
		15 & 12 & [0.2298]  &  174       &  89        & 98 \\
		15 & 14 & [0.6535]  &  [0.0048]  &  [0.0032]  & [0.0077] \\
		15 & 15 & [1.4129]  &  [0.0073]  &  [0.0084]  & [0.0102] \\
		15 & 16 & [1.1820]  &  [0.0125]  &  [0.0121]  & [0.0136] \\
		15 & 17 & [4.0515]  &  [0.0119]  &  [0.0115]  & [0.0139] \\
		15 & 18 & [6.4008]  &  [0.0115]  &  [0.0096]  & [0.0096] \\
		15 & 20 & [8.8120]  &  [0.0111]  &  [0.0158]  & [0.0095] \\
		\hline
		& Average & 2.8550 & 0.0074 & 0.0076 & 0.0081  \\
		\hline
	\end{tabular}}
	\label{table-2class}
\end{table}

The results displayed in Table \ref{table-2class} show unequivocally that the three cutting plane approaches {\tt SCA}, {\tt B-SCA}, and {\tt OA-B-SCA} dominate a direct solution of {\bf MSP}. For instances with three searchers, the average optimality gap of each cutting plane method is below 10\% while the one obtained by solving directly {\bf MSP} exceeds 50\%.
Comparing the cutting plane algorithms, we see that {\tt SCA}, {\tt B-SCA}, and {\tt OA-B-SCA} exhibit similar performance levels for the relatively easy instances (i.e. $J=15$). For challenging cases involving $J=3$ searchers, {\tt OA-B-SCA} performs much better, on average {\tt SCA} and {\tt B-SCA} produce twice as large optimality gaps. The results in Table \ref{table-2class} demonstrate that the proposed {\tt OA-B-SCA} is most effective for the most challenging instances.

Next, we consider Table \ref{table-2class3} where the searchers vary in both endurance and detection ability, and thus the rate modification factors $\beta_{l,s',s,t}$ cannot all be 1. The detection ability of class-two searchers is equal to 80\% of that of class-one searchers. The resulting instances are exceptionally challenging, in particular when the numbers of periods and searchers increase. The cutting plane method  {\tt OA-B-SCA} is the most efficient approach as it provides by far the smallest optimality gap for each instance, and is the only method that can solve one instance to optimality within 900 seconds. It provides practically reasonable optimality gaps for planning horizon $T\leq 12$. Analysis of each instance reveals that the high optimality gap for {\bf MSP} is usually due to the weakness of its lower bound. For example, the best lower bound  for the $J=15$, $T=12$ instance -- obtained by {\tt OA-B-SCA} -- confirms that the best integer solution (i.e., with objective value of 0.3455) from {\bf MSP} actually has an optimality gap of 7\%. This is dramatically better than the 98\% reported in Table \ref{table-2class3}. Thus, {\bf MSP} cannot be ruled out as a viable approach for generating good feasible solutions.

\begin{table}[H]
	\centering
	\caption{{\small For Markovian target model, $S=83$ states, and varying numbers of time periods and searchers across two classes with varying endurance and detection ability: Upper bound (UB) and solution time (sec.) to relative optimality gap of 0.0001 or, if not reached in 900 seconds, relative optimality gap in brackets after 900 seconds ($\infty$ indicates that no bound is available).}}
	{\small\begin{tabular}{|c|c|r|r|r|r|r|r|}
		\hline
		&  & \multicolumn{2}{c|}{{\bf MSP}} & \multicolumn{2}{c|}{{\tt B-SCA}} & \multicolumn{2}{c|}{{\tt OA-B-SCA}}\\
		\hline
		$J$ & $T$ & Time     & UB     & Time     & UB     & Time   & UB\\
		\hline
		3 & 10  & [0.1552] & 0.4245 & [0.1195]  & 0.4369  & [0.0001] & 0.3779\\
		3 & 12  & [0.4067] & 0.3742 & [0.6347]  & 0.4534  & [0.0024] & 0.3230\\
		3 & 14  & [0.8468] & 0.3178 & [3.3363]  & 0.4380  & [0.0913] & 0.2331\\
		3 & 16  & [1.0000] & 0.5263 & [115.65]  & 0.4229  & [0.1983] & 0.2094\\
		3 & 18  & [1.2159] & 0.7555 & [$\infty$]& 0.3557  & [0.2767] & 0.1649\\
		\hline
		15 & 10 & [0.5524] & 0.3900 & [0.0286]  & 0.3885  & 112      & 0.3791\\
		15 & 12 & [0.9791] & 0.3455 & [0.0726]  & 0.3336  & [0.0272] & 0.3314 \\
		15 & 14 & [1.4006] & 0.6878 & [0.5835]  & 0.2639  & [0.3494] & 0.2778 \\
		15 & 16 & [1.4425] & 0.8953 & [1.9474]  & 0.2673  & [1.0097] & 0.2727 \\
		15 & 18 & [1.9703] & 0.8542 & [269.44]  & 0.2236  & [1.9238] & 0.2378 \\
		\hline
	\end{tabular}}
	\label{table-2class3}
\end{table}

\section{Conclusion}

Search planning for a randomly moving target in discrete time and space should account for operationally important concerns such as the employment of heterogeneous searchers with distinct endurance level, detection ability, and travel speed, the need for deconfliction among the searchers, and the ability for the target to camouflage and thus making any sensor ineffective. We account for all these concerns within a convex mixed-integer nonlinear program, while taking advantage of homogeneous sensors and Markovian target path models when present.

Since the objective function is a weighted sum of exponential functions with integer arguments, it can be linearized. We propose a new linearization technique and extend two existing ones to account for heterogeneous searchers and operational constraints. While equivalent to the actual problem, the linearizations tend to be large-scaled but reducible via customized preprocessing and lazy-constraint techniques. We also develop three cutting plane methods for challenging instances. The most suitable approach for a particular problem instance depends on the number of searchers, the length of the planning horizon, and, maybe primarily, on the characteristics of the target movement.

When the target follows any one of a moderately large number of paths (e.g., 1000 paths), it turns out that a direct solution of a linearization (after preprocessing) by a standard mixed-integer linear programming solver is viable and in fact computationally most effective as long as the searchers are essentially homogeneous and the planning horizon is no longer than 15 time periods. For example, an instance with 82 states, 15 time periods, 50 homogeneous searchers, no endurance constraints, and no camouflaging solves to optimality in less than one minute using Gurobi. For more complex instances involving heterogeneous searchers, our lazy-constraint-based outer-approximation algorithm becomes the most efficient approach.
When the target moves according to a Markov chain, which tends to produce a massive number of possible paths, the linearizations become inefficient and we rely on three cutting plane methods. Two of these are complemented with a bundle approach and the last one is embedded in an outer-approximation algorithm. The latter performs best on instances with heterogeneous searchers. For example, we achieve an optimality gap of 2.7\% after 900 seconds for an instance with 83 states, 12 time periods, a camouflaging target, and 15 searchers across two classes of varying sensor capabilities and endurance.

Our extensive numerical study also provides some insights for practitioners regarding the impact of endurance, detection ability, and camouflage. Searchers facing endurance limitations tend to delay the search and wait for the target to approach them to avoid wasting time in transit to the target's likely location. Increased travel speed for the searchers improves the probability of detecting the target, but possibly only with a moderate amount. A camouflaging target is less mobile and results in a concentrated search plan near the target's initial position.

\section*{Acknowledgement}
Lejeune acknowledges support from NSF (ECCS-2114100 and RISE-2220626) and ONR (N00014-22-1-2649); Royset acknowledges support from ONR (N0001423WX01316; N0001423WX00403).

\bibliographystyle{plain}
\bibliography{refs}

\end{document}